# ESTIMATING LINEAR FUNCTIONALS IN NONLINEAR REGRESSION WITH RESPONSES MISSING AT RANDOM

By Ursula U. Müller

*Texas A&M University*

We consider regression models with parametric (linear or nonlinear) regression function and allow responses to be "missing at random." We assume that the errors have mean zero and are independent of the covariates. In order to estimate expectations of functions of covariate and response we use a *fully imputed* estimator, namely an empirical estimator based on estimators of conditional expectations given the covariate. We exploit the independence of covariates and errors by writing the conditional expectations as unconditional expectations, which can now be estimated by empirical plug-in estimators. The mean zero constraint on the error distribution is exploited by adding suitable residual-based weights. We prove that the estimator is efficient (in the sense of Hájek and Le Cam) if an efficient estimator of the parameter is used. Our results give rise to new efficient estimators of smooth transformations of expectations. Estimation of the mean response is discussed as a special (degenerate) case.

**1. Introduction.** Consider a regression model $Y = r_\vartheta(X) + \varepsilon$ with linear or nonlinear regression function $r_\vartheta$ depending on a finite-dimensional parameter $\vartheta$ in some open set. Assume that the covariate vector $X$ and the error variable $\varepsilon$ are independent and that $E\varepsilon = 0$. Note that we do not make any further model assumptions on the distributions of the variables. We are interested in the situation where the response $Y$ is *missing at random*, in other words, we always observe $X$ but only observe $Y$ in those cases where some indicator $Z$ equals one, and the indicator $Z$ is conditionally independent of $Y$ given $X$.

We want to estimate the expectation $Eh(X, Y)$ of some known square-integrable function $h$ from a sample $(X_i, Z_i Y_i, Z_i)$, $i = 1, \ldots, n$, for example, the mean response, higher moments of $Y$ or $X$ or mixed moments. If all indicators $Z_i$ were 1, a simple consistent estimator would be the empirical estimator $n^{-1}\sum_{i=1}^{n} h(X_i, Y_i)$. A related estimator for the missing data situation considered here would be

$$\frac{1}{n}\sum_{i=1}^{n}\frac{Z_i}{\hat{\pi}(X_i)}h(X_i, Y_i)$$

with $\hat{\pi}(X)$ denoting an estimator of the conditional probability $\pi(X) = P(Z = 1|X) = E(Z|X)$. Another estimator is the *partially imputed* estimator

$$\frac{1}{n}\sum_{i=1}^{n}\{Z_i h(X_i, Y_i) + (1 - Z_i)\hat{\chi}(X_i)\},$$

where $\hat{\chi}(X)$ is a (semiparametric) estimator of the conditional expectation $\chi(X) = E\{h(X, Y)|X\}$. An alternative to this estimator is the *fully imputed* estimator $n^{-1}\sum_{i=1}^{n}\hat{\chi}(X_i)$.

If a nonparametric estimator $\hat{\chi}$ is used, we expect all three estimators to be asymptotically equivalent. For $h(X, Y) = Y$ and the last two estimators, this is sketched in Cheng (1994). Here we assume a specific form of the conditional distribution of $Y$ given $X$, and we can construct better estimators than the nonparametric ones. We then expect the *fully imputed* estimator $n^{-1}\sum_{i=1}^{n}\hat{\chi}(X_i)$ to be better than the partially imputed one, which in turn should be better than the first estimator. For *parametric* models this is shown for $h(X, Y) = Y$ by Tamhane (1978) and Matloff (1981). Müller, Schick and Wefelmeyer (2006) show for several regression models (not including the present one) and arbitrary $h$ that the fully imputed estimator is usually better than the partially imputed estimator. That the same holds for the nonlinear regression model considered here is intuitively clear: our model $E(Y|X) = r_\vartheta(X)$ constitutes a structural constraint. The fully imputed estimator, based on estimators $\hat{\chi}(X)$ that use the structure, will therefore be better than the partially imputed estimator, which uses this information only at data points where responses are missing.

In this article we study the fully imputed estimator based on suitable estimators for $\chi(X)$ and show that it is efficient. The construction is as follows: in a first step we exploit the independence of covariates and errors and the structure of the regression model and write the conditional expectation $\chi(x) = \chi(x, \vartheta)$ as an unconditional expectation of the error distribution,

$$\chi(x, \vartheta) = E\{h(X, Y)|X = x\}$$
$$= Eh\{x, r_\vartheta(x) + \varepsilon\} = Eh\{x, r_\vartheta(x) + Y - r_\vartheta(X)\}.$$

This representation suggests an empirical plug-in estimator based on the observed data, namely

$$\hat{\chi}(x, \hat{\vartheta}) = \sum_{j=1}^{n} Z_j h\{x, r_{\hat{\vartheta}}(x) + Y_j - r_{\hat{\vartheta}}(X_j)\}\Big/\sum_{j=1}^{n} Z_j,$$

where $\hat{\vartheta}$ is an estimator of $\vartheta$. The corresponding fully imputed estimator is

$$(1.1) \qquad \frac{1}{n}\sum_{i=1}^{n}\hat{\chi}(X_i,\hat{\vartheta}) = \frac{1}{n}\sum_{i=1}^{n}\frac{\sum_{j=1}^{n}Z_j h\{X_i, r_{\hat{\vartheta}}(X_i) + Y_j - r_{\hat{\vartheta}}(X_j)\}}{\sum_{j=1}^{n}Z_j}.$$

It is straightforward to check that $\hat{\chi}(x,\vartheta)$ is consistent for $Eh\{x, r_{\vartheta}(x) + \varepsilon\}$ [which yields consistency of $n^{-1}\sum_{i=1}^{n}\hat{\chi}(X_i,\hat{\vartheta})$, with $\hat{\vartheta}$ consistent]; note that $\hat{\chi}(x,\vartheta)$ tends in probability to $E[Zh\{x, r_{\vartheta}(x) + \varepsilon\}]/EZ$ with $EZ = E\{E(Z|X)\} = E\pi(X)$. Now use the *missing at random* assumption and the independence of $X$ and $\varepsilon$ to rewrite the numerator,

$$E(E[Zh\{x, r_{\vartheta}(x) + \varepsilon\}|X]) = E(E(Z|X)E[h\{x, r_{\vartheta}(x) + \varepsilon\}|X])$$
$$= E[\pi(X)Eh\{x, r_{\vartheta}(x) + \varepsilon\}]$$
$$= E\pi(X)Eh\{x, r_{\vartheta}(x) + \varepsilon\}.$$

The limit of $\hat{\chi}(x,\vartheta)$ is therefore $\chi(x,\vartheta) = Eh\{x, r_{\vartheta}(x) + \varepsilon\}$.

The estimator (1.1) is well thought out and consistent. However, it is not yet efficient, even if an efficient estimator for $\vartheta$ is used (which is relatively elaborate in the model considered here; see Section 5): we focus on the common situation where the errors have mean zero; this information must also be incorporated in order to obtain efficiency.

Motivated by Owen's empirical likelihood approach, we improve the above estimator by introducing weights which use the mean zero constraint on the error distribution. However, and in contrast to the original approach, we cannot observe the errors and must use residuals. This clearly complicates the situation: since we have missing responses the residuals are partially incomplete and, moreover, they involve parameter estimates $\hat{\vartheta}$. Formally, we choose weights $\hat{w}_j$ based on residuals $\hat{\varepsilon}_j = Y_j - r_{\hat{\vartheta}}(X_j)$ such that $\sum_{j=1}^{n}\hat{w}_j Z_j \hat{\varepsilon}_j = 0$. (See Section 3 for more details.)

Our final estimator now is a weighted version of the above fully imputed estimator, namely

$$(1.2) \qquad \frac{1}{n}\sum_{i=1}^{n}\hat{\chi}_w(X_i,\hat{\vartheta}) = \frac{1}{n}\sum_{i=1}^{n}\frac{\sum_{j=1}^{n}\hat{w}_j Z_j h\{X_i, r_{\hat{\vartheta}}(X_i) + Y_j - r_{\hat{\vartheta}}(X_j)\}}{\sum_{j=1}^{n}Z_j}.$$

The combination of full imputation methods (involving estimators of *unconditional* expectations of the error distribution) with empirical likelihood ideas provides a new methodology which has not appeared in the literature before. We show in this article that $n^{-1}\sum_{i=1}^{n}\hat{\chi}_w(X_i,\hat{\vartheta})$ is efficient if an efficient estimator $\hat{\vartheta}$ for $\vartheta$ is used. The partially imputed estimator will in general not be efficient, even if $\hat{\vartheta}$ is efficient for $\vartheta$.

For estimation of the mean response, that is, if $h(X,Y) = Y$, which is of particular interest and typically considered in the literature, the estimator simplifies to the straightforward estimator $n^{-1}\sum_{i=1}^{n}r_{\hat{\vartheta}}(X_i)$. That the unweighted estimator (1.1) for $EY$ cannot be efficient is immediately apparent: consider the case where all responses are observed. Here (1.1) reduces to the empirical estimator $n^{-1} \sum_{i=1}^{n} Y_i$ which does not use the regression structure at all. It will be seen that its influence function is not the efficient one. (See Section 6 for details.)

Our efficiency results are based on the Hájek–Le Cam theory for locally asymptotically normal families. As a consequence, our proposed estimators have a limiting normal distribution with the asymptotic variance determined by the influence function. It is therefore straightforward to construct asymptotic confidence interval for $Eh(X, Y)$ (see Section 6.3).

In addition, estimators for smooth (continuously differentiable) transformations of expectations $Eh(X, Y)$ are also now available, with the variance of the response, $\text{Var}\, Y = EY^2 - E^2Y$, as an important example. Since efficiency is preserved by smooth transformations, plugging in efficient estimators yields an efficient estimator of the transformation. The transformation for $\text{Var}\, Y$ in terms of the first two moments is $(EY, EY^2) \mapsto EY^2 - (EY)^2$. Plugging in $n^{-1} \sum_{i=1}^{n} r_{\hat{\vartheta}}(X_i)$ for $EY$ and the weighted fully imputed estimator for $EY^2$ (which is straightforward to compute and is also given in Section 6) gives an efficient estimator of the variance.

To our knowledge, our estimator (1.2) is the first efficient estimator for arbitrary linear functionals $Eh(X, Y)$ (including the mean functional $EY$) in the nonlinear regression model (including the linear regression model $Y = \vartheta^{\top} X + \varepsilon$) with independent centered errors when responses are missing at random. Matloff (1981) considers estimation of the mean $EY$ in a model related to ours, the (parametric) *conditional mean model*, $E(Y|X) = r_{\vartheta}(X)$, which can (but need not) also be written in the form $Y = r_{\vartheta}(X) + \varepsilon$ with conditionally centered errors, $E(\varepsilon|X) = 0$. He shows that the average of the estimated regression function values (with his estimator $\hat{\vartheta}$ of $\vartheta$) improves upon the partially imputed estimator. Wang and Rao (2001) consider linearly constrained covariates and develop an empirical likelihood approach for inference about the mean in linear regression (with independent errors) based on partial linear regression imputation. In Wang and Rao (2002) they present an empirical likelihood approach for inference about the mean response in nonparametric regression, based on partial kernel regression imputation as suggested by Cheng (1994). A different empirical likelihood method for this setting is proposed by Qin and Zhang (2007). Wang (2004) assumes a parametric model for the conditional density of $Y$ given $X$, with constraints on the covariate distribution, and introduces a weighted partial imputation estimator for the mean, utilizing empirical likelihood techniques. Wang, Linton and Härdle (2004) consider a partially linear regression model for the conditional mean function and derive inference tools for the mean response based on a class of asymptotically equivalent (partially and fully imputed) estimators. A related article is Liang, Wang and Carroll (2007) who additionally assume that covariates are measured with error. Chen, Fan, Li and Zhou (2006) consider partially imputed estimators for the mean response in a quasi-likelihood setting. Maity, Ma and Carroll (2007) estimate expectations in semi-parametric regression models, with and without missing responses. They consider a general regression function involving a parametric and a nonparametric part, thus covering the partly linear model, and assume that the likelihood function given the covariates is known.

For estimating expectations, little attention has been given to the fully imputed estimator. We anticipate that in many situations, in particular in models with structural assumptions, improved estimators can be obtained by using appropriate full imputation instead of partial imputation estimates.

Inference for missing data has been studied by many authors, also recently. Chen and Wang (2009) study estimation of parameters which are defined by model constraints. They introduce an empirical likelihood approach involving estimating equations, where missing variables are replaced using a nonparametric imputation approach. Chen, Hong and Tarozzi (2008) consider parameter estimation as well. They introduce efficient estimators for parameters in GMM models with missing data, and assume that the missingness can be explained by auxiliary variables. More references to recent literature can be found, for example, in Wang, Linton and Härdle (2004) and in the monograph by Tsiatis (2006). For an introduction, see Tsiatis (2006) and the books by Little and Rubin (2002) and Gelman et al. (1995).

This paper is organized as follows. In Section 2 we derive a stochastic expansion of the unweighted estimator. The expansion of the weighted estimator is given in Section 3, utilizing the results of Section 2. Section 4 characterizes efficient estimators of arbitrary functionals of the joint distribution and gives the efficient influence function of the functional $Eh(X, Y)$ in the nonlinear regression model. In Section 5 we characterize efficient estimators for the parameter vector $\vartheta$ and briefly sketch the construction of such an estimator. In this section we also show our main result, that the weighted estimator with an efficient estimator $\hat{\vartheta}$ for $\vartheta$ plugged in is efficient for $Eh(X, Y)$. Section 6 contains a short discussion of special cases such as estimation of the mean response. We also compare, using computer simulations, the efficient (weighted fully imputed) estimator with the other approaches, with convincing results. For these studies we considered a linear and a nonlinear regression function and estimation of two simple functionals, namely of the response mean and second moment, for which the efficient (weighted fully imputed) estimator simplifies, and estimation of a more complicated expectation. We also briefly sketch the construction of confidence intervals.

**2. Expansion of the unweighted estimator.** In this section we derive an expansion of the unweighted estimator $n^{-1} \sum_{i=1}^{n} \hat{\chi}(X_i, \hat{\vartheta})$, which is a special case of the weighted estimator $n^{-1} \sum_{i=1}^{n} \hat{\chi}_w(X_i, \hat{\vartheta})$ with all weights being equal to one, $w_j = 1$. This can be regarded as a result of independent interest since the estimator (with an appropriate estimator $\hat{\vartheta}$) would be relevant for regression models where the errors cannot be assumed to have mean zero. Also, we will see in the next section that the weighted estimator can be written as the sum of the unweighted estimator and an additional correction term. Hence we can utilize the results later when we derive an expansion of the weighted estimator.

Throughout this paper we will assume that $Y$ is square integrable and that the error variance $E\varepsilon^2 = \sigma^2$ is nonzero and finite. We also suppose that the error distribution has a Lebesgue density $f$ and finite Fisher information, $E\ell^2(\varepsilon) < \infty$, where $\ell$ denotes the score function for location, $\ell(\varepsilon) = -f'(\varepsilon)/f(\varepsilon)$. The degenerate case that we (almost surely) never observe a response $Y$ will be excluded by assuming $P(Z=1) = EZ > 0$. The following assumptions will also be required.

ASSUMPTION 1. The regression function $\tau \mapsto r_\tau(x)$ is differentiable at $\tau = \vartheta$ with a $p$-dimensional square integrable gradient $\dot{r}_\vartheta(x)$ which satisfies the Lipschitz condition

$$|\dot{r}_\tau(x) - \dot{r}_\vartheta(x)| \le |\tau - \vartheta| a(x), \qquad a(X) \text{ square integrable.}$$

Later we will also need that the covariance matrix of an efficient parameter estimator $\hat{\vartheta}$ [which involves the covariance matrix of $\dot{r}_\vartheta(X)$ and the Fisher information] is invertible.

Now use a Taylor expansion to see that

$$\sum_{i=1}^{n} \{r_\tau(X_i) - r_\vartheta(X_i) - \dot{r}_\vartheta(X_i)^\top (\tau - \vartheta)\}^2$$

$$= \sum_{i=1}^{n} \left[ \int_0^1 \{\dot{r}_{\vartheta + u(\tau - \vartheta)}(X_i) - \dot{r}_\vartheta(X_i)\}^\top (\tau - \vartheta) \, du \right]^2$$

$$\le |\tau - \vartheta|^2 \sum_{i=1}^{n} \int_0^1 |\dot{r}_{\vartheta + u(\tau - \vartheta)}(X_i) - \dot{r}_\vartheta(X_i)|^2 \, du$$

$$\le |\tau - \vartheta|^4 \sum_{i=1}^{n} a^2(X_i).$$

Assumption 1 therefore guarantees that the function $\tau \mapsto r_\tau(X)$ is stochastically differentiable, that is, for each constant $C$,

$$(2.1) \qquad \sup_{|\tau - \vartheta| \le C n^{-1/2}} \sum_{i=1}^{n} \{r_\tau(X_i) - r_\vartheta(X_i) - \dot{r}_\vartheta(X_i)^\top (\tau - \vartheta)\}^2 = o_p(1).$$

We will not need the first partial derivative of $h(x,y)$, $\partial/\partial x h(x,y)$. Therefore we will write $h'$ for the second partial derivative, $h'(x,y) = \partial_2 h(x,y) = \partial/\partial y h(x,y)$.

ASSUMPTION 2. *The function $h(x,y)$ is differentiable in $y$ with a square integrable partial derivative $h'(x,y) = \partial/\partial y h(x,y)$ which satisfies the Lipschitz condition*

$$|h'(x,z) - h'(x,y)| \le |z-y| b(x,y), \qquad b(X,Y) \text{ square integrable}.$$

In the following $\bar{Z}$ will denote the average of the indicators $Z_i$, $\bar{Z} = n^{-1} \sum_{i=1}^{n} Z_i$. The next lemma gives the expansion of the estimator around the true parameter $\vartheta$.

LEMMA 2.1. *Assume that Assumptions 1 and 2 hold and that $\hat{\vartheta}$ is a $\sqrt{n}$ consistent estimator of $\vartheta$. Then the unweighted estimator has the expansion*

$$(2.2) \qquad \frac{1}{n}\sum_{i=1}^{n}\hat{\chi}(X_i,\hat{\vartheta}) = \frac{1}{n}\sum_{i=1}^{n}\hat{\chi}(X_i,\vartheta) + D^{\top}(\hat{\vartheta} - \vartheta) + o_p(n^{-1/2})$$

*with $D = E(h(X,Y)[\dot{r}_\vartheta(X) - E\{\dot{r}_\vartheta(X)|Z=1\}]\ell(\varepsilon))$.*

PROOF. For reasons of clarity we introduce the notation

$$f_{ij}(\vartheta) = h\{X_i, r_\vartheta(X_i) + Y_j - r_\vartheta(X_j)\}$$

and write $\dot{f}_{ij}$ for the gradient. Then

$$\frac{1}{n}\sum_{i=1}^{n}\hat{\chi}(X_i,\hat{\vartheta})$$

$$= \frac{1}{n^2}\sum_{i=1}^{n}\sum_{j=1}^{n}\frac{Z_j}{\bar{Z}}h\{X_i, r_{\hat{\vartheta}}(X_i) + Y_j - r_{\hat{\vartheta}}(X_j)\}$$

$$= \frac{1}{\bar{Z}}\frac{1}{n^2}\sum_{i=1}^{n}\left[\sum_{\substack{j=1\\j\ne i}}^{n}Z_j h\{X_i, r_{\hat{\vartheta}}(X_i) + Y_j - r_{\hat{\vartheta}}(X_j)\} + Z_i h(X_i,Y_i)\right]$$

$$(2.3)$$

$$= \frac{1}{\bar{Z}}\frac{1}{n^2}\sum_{i=1}^{n}\left\{\sum_{\substack{j=1\\j\ne i}}^{n}Z_j f_{ij}(\vartheta) + Z_i h(X_i,Y_i)\right\}$$

$$\qquad + \frac{1}{\bar{Z}}\frac{1}{n^2}\sum_{i=1}^{n}\sum_{\substack{j=1\\j\ne i}}^{n}Z_j\{f_{ij}(\hat{\vartheta}) - f_{ij}(\vartheta)\}$$

$$= \frac{1}{n}\sum_{i=1}^{n}\hat{\chi}(X_i,\vartheta) + \frac{1}{\bar{Z}}\frac{1}{n^2}\sum_{i=1}^{n}\sum_{\substack{j=1\\j\ne i}}^{n}Z_j\{f_{ij}(\hat{\vartheta}) - f_{ij}(\vartheta)\}.$$

Below we will show that

$$(2.4) \qquad \frac{1}{Z} \frac{1}{n^2} \sum_{i=1}^{n} \sum_{\substack{j=1 \\ j \neq i}}^{n} Z_j \{f_{ij}(\hat{\vartheta}) - f_{ij}(\vartheta)\} = D^\top(\hat{\vartheta} - \vartheta) + o_p(n^{-1/2})$$

with $D = (EZ)^{-1} E[Z_2 h'\{X_1, r_\vartheta(X_1) + Y_2 - r_\vartheta(X_2)\} \{\dot{r}_\vartheta(X_1) - \dot{r}_\vartheta(X_2)\}]$. That this $D$ is indeed of the form given in the lemma can be seen as follows. Consider

$$D = E[h'\{X_1, r_\vartheta(X_1) + \varepsilon_2\} \dot{r}_\vartheta(X_1)]$$
$$- \frac{1}{EZ} E[h'\{X_1, r_\vartheta(X_1) + \varepsilon_2\} \mathbf{1}(Z_2 = 1) \dot{r}_\vartheta(X_2)].$$

The first term can be written $E(E[h'\{X_1, r_\vartheta(X_1) + \varepsilon_2\}|X_1] \dot{r}_\vartheta(X_1))$. Integration by parts of the inner integral gives $E[h'\{X_1, r_\vartheta(X_1) + \varepsilon_2\}|X_1] = E[h\{X_1, r_\vartheta(X_1) + \varepsilon_2\} \ell(\varepsilon_2)|X_1]$. The second term is $E[h'\{X_1, r_\vartheta(X_1) + \varepsilon_2\}] E\{\dot{r}_\vartheta(X)|Z = 1\}$. We proceed analogously and, in conclusion, obtain

$$(2.5) \qquad D = E(h(X, Y)[\dot{r}_\vartheta(X) - E\{\dot{r}_\vartheta(X)|Z = 1\}] \ell(\varepsilon)).$$

The result now follows from (2.3), (2.4) and (2.5). It remains to verify (2.4). The proof consists of two parts,

$$(2.6) \qquad \frac{1}{Z} \frac{1}{n^2} \sum_{i=1}^{n} \sum_{\substack{j=1 \\ j \neq i}}^{n} Z_j \{f_{ij}(\hat{\vartheta}) - f_{ij}(\vartheta) - \dot{f}_{ij}(\vartheta)^\top(\hat{\vartheta} - \vartheta)\} = o_p(n^{-1/2}),$$

$$(2.7) \qquad \frac{1}{Z} \frac{1}{n^2} \sum_{i=1}^{n} \sum_{\substack{j=1 \\ j \neq i}}^{n} Z_j \dot{f}_{ij}(\vartheta)^\top(\hat{\vartheta} - \vartheta) = D^\top(\hat{\vartheta} - \vartheta) + o_p(n^{-1/2}).$$

Statement (2.7) can be quickly proved: since $\hat{\vartheta}$ is $\sqrt{n}$ consistent we can replace the gradient by its expectation,

$$\frac{1}{Z} \frac{1}{n^2} \sum_{i=1}^{n} \sum_{\substack{j=1 \\ j \neq i}}^{n} Z_j \dot{f}_{ij}(\vartheta)^\top(\hat{\vartheta} - \vartheta)$$

$$= \frac{1}{Z} \frac{1}{n^2} \sum_{i=1}^{n} \sum_{\substack{j=1 \\ j \neq i}}^{n} E\{Z_j \dot{f}_{ij}(\vartheta)\}^\top(\hat{\vartheta} - \vartheta) + o_p(n^{-1/2})$$

$$= \frac{1}{EZ} E\{Z_2 \dot{f}_{12}(\vartheta)\}^\top(\hat{\vartheta} - \vartheta) + o_p(n^{-1/2})$$

with $(EZ)^{-1} E\{Z_2 \dot{f}_{12}(\vartheta)\} = D^\top$ as given in (2.4). For the proof of (2.6) it suffices to show that

$$\sum_{i=1}^{n} \left[ \frac{1}{\sqrt{n}} \sum_{\substack{j=1 \\ j \neq i}}^{n} Z_j \{f_{ij}(\hat{\vartheta}) - f_{ij}(\vartheta) - \dot{f}_{ij}(\vartheta)^\top(\hat{\vartheta} - \vartheta)\} \right]^2 = O_p(1).$$

This holds by the following arguments. Rewrite the above expression and apply the Cauchy–Schwarz inequality to obtain

$$\sum_{i=1}^{n}\left(\frac{1}{\sqrt{n}}\sum_{\substack{j=1\\j\neq i}}^{n}Z_j\int_0^1[\dot{f}_{ij}\{\vartheta+u(\hat{\vartheta}-\vartheta)\}-\dot{f}_{ij}(\vartheta)]^\top(\hat{\vartheta}-\vartheta)\,du\right)^2$$

$$\leq\sum_{i=1}^{n}\sum_{\substack{j=1\\j\neq i}}^{n}Z_j|\hat{\vartheta}-\vartheta|^2\int_0^1|\dot{f}_{ij}\{\vartheta+u(\hat{\vartheta}-\vartheta)\}-\dot{f}_{ij}(\vartheta)|^2\,du.$$

The difference $|\dot{f}_{ij}\{\vartheta+u(\hat{\vartheta}-\vartheta)\}-\dot{f}_{ij}(\vartheta)|^2$ is bounded by $|\hat{\vartheta}-\vartheta|^2$ times a square integrable function $A_{ij}$. This holds due to Assumptions 1 and 2, namely the Lipschitz conditions on $\dot{r}_\vartheta$ and $h'$ and since $a(X), b(X,Y), \dot{r}_\vartheta(X)$ and $h'(X,Y)$ are square integrable. Summing up, the expression is bounded by $|\hat{\vartheta}-\vartheta|^4\sum_{i=1}^{n}\sum_{j=1,j\neq i}^{n}A_{ij}$ which is stochastically bounded since $\hat{\vartheta}$ is $\sqrt{n}$ consistent. $\square$

We will now replace the estimated conditional expectation $\hat{\chi}$ in the right-hand side of (2.2) by the true one. Set

$$S=\frac{1}{n(n-1)}\sum_{i=1}^{n}\sum_{\substack{j=1\\j\neq i}}^{n}\frac{Z_j}{EZ}h\{X_i,r_\vartheta(X_i)+Y_j-r_\vartheta(X_j)\}.$$

We have

$$\frac{1}{n}\sum_{i=1}^{n}\hat{\chi}(X_i,\vartheta)=\frac{EZ}{\bar{Z}}S+O_p(n^{-1})=S-\frac{\bar{Z}-EZ}{EZ}ES+o_p(n^{-1/2})$$

and, by the Hoeffding decomposition,

$$S=ES+\frac{1}{n}\sum_{i=1}^{n}\{\chi(X_i,\vartheta)-ES\}+\frac{1}{n}\sum_{i=1}^{n}\left\{\frac{Z_i\bar{h}(\varepsilon_i)}{EZ}-ES\right\}+o_p(n^{-1/2})$$

with $\bar{h}(\varepsilon)=E\{h(X,Y)|\varepsilon\}$, $ES=Eh(X,Y)=E\bar{h}(\varepsilon)$. Combining the above yields

$$\frac{1}{n}\sum_{i=1}^{n}\hat{\chi}(X_i,\vartheta)=\frac{1}{n}\sum_{i=1}^{n}\chi(X_i,\vartheta)+\frac{1}{n}\sum_{i=1}^{n}\frac{Z_i}{EZ}\{\bar{h}(\varepsilon_i)-E\bar{h}(\varepsilon)\}+o_p(n^{-1/2}).$$

This and Lemma 2.1 give our expansion for the unweighted estimator which we formulate as a corollary.

COROLLARY 2.2. *Assume that Assumptions 1 and 2 hold and that $\hat{\vartheta}$ is a $\sqrt{n}$ consistent estimator of $\vartheta$. Then, with $D=E(h(X,Y)[\dot{r}_\vartheta(X)-$*

$E\{\dot{r}_\vartheta(X)|Z=1\}\ell(\varepsilon))$ *and* $\bar{h}(\varepsilon) = E\{h(X,Y)|\varepsilon\}$, *the unweighted estimator has the expansion*

$$\frac{1}{n}\sum_{i=1}^n \hat{\chi}(X_i, \hat{\vartheta})$$
$$= \frac{1}{n}\sum_{i=1}^n \left[\chi(X_i, \vartheta) + \frac{Z_i}{EZ}\{\bar{h}(\varepsilon_i) - E\bar{h}(\varepsilon_i)\}\right] + D^\top(\hat{\vartheta} - \vartheta) + o_p(n^{-1/2}).$$

**3. Expansion of the weighted estimator.** In this section we study the weighted estimator which uses residual-based weights, $\hat{w}_j$, that are constructed by adapting empirical likelihood techniques. The approach is to maximize $\prod_{j=1}^n \hat{w}_j$ subject to the mean zero constraint on the error distribution, $\sum_{j=1}^n \hat{w}_j Z_j \hat{\varepsilon}_j = 0$, with $\hat{w}_j \geq 0$ and $\sum_{j=1}^n \hat{w}_j = n$. The weights solving this optimization problem are given by $\hat{w}_j = 1/(1 + \hat{\lambda} Z_j \hat{\varepsilon}_j)$, where $\hat{\lambda}$ denotes the Lagrange multiplier—provided $\hat{\lambda}$ exists. As shown by Owen ([1988](#), [2001](#)), this is the case if not all residuals have the same sign, that is, on the event $\min_{1 \leq j \leq n} \hat{\varepsilon}_j < 0 < \max_{1 \leq j \leq n} \hat{\varepsilon}_j$, which has probability tending to one since the residuals $\hat{\varepsilon}_j$ are uniformly close to the centered errors $\varepsilon_j$ [see [(A.1)](#) in the [Appendix](#)]. If $\hat{\lambda}$ does not exist, we set $\hat{\lambda} = 0$. Note that the weights equal one if $Z_j = 0$ or $\hat{\lambda} = 0$. For computational issues we refer to Section 2.9 of Owen's book ([2001](#)).

The formula for the weights can be written as an identity, $\hat{w}_j = 1 - \hat{\lambda}\hat{w}_j Z_j \hat{\varepsilon}_j$. This enables us to decompose the estimator into the unweighted estimator and an additional correction term,

$$(3.1) \quad \begin{aligned} \frac{1}{n}\sum_{i=1}^n \hat{\chi}_w(X_i, \hat{\vartheta}) &= \frac{1}{n}\sum_{i=1}^n \hat{\chi}(X_i, \hat{\vartheta}) \\ &\quad - \frac{\hat{\lambda}}{n^2}\sum_{i=1}^n \sum_{j=1}^n \hat{w}_j \hat{\varepsilon}_j \frac{Z_j}{Z} h\{X_i, r_{\hat{\vartheta}}(X_i) + \hat{\varepsilon}_j\}. \end{aligned}$$

Since we have already derived an expansion of the unweighted estimator (see Corollary [2.2](#)) we only need to study the second term on the right-hand side. In Lemma [3.1](#) we will derive an expansion of the estimated Lagrange multiplier $\hat{\lambda}$ and use this result in Lemma [3.2](#), where we determine an approximation of the extra term. For the proof of Lemma [3.1](#) we proceed analogously to Owen (2001), pages 219–221 [compare also Müller, Schick and Wefelmeyer [(2005)](#)]. This requires some auxiliary results which are proved in the [Appendix](#), namely

$$(3.2) \quad \max_{1 \leq i \leq n} |Z_i \hat{\varepsilon}_i| = o_p(n^{1/2}),$$

$$
(3.3) \quad \frac{1}{n} \sum_{i=1}^{n} Z_i \hat{\varepsilon}_i = \frac{1}{n} \sum_{i=1}^{n} Z_i \varepsilon_i - E Z E \{\dot{r}_\vartheta(X) | Z = 1\}^\top (\hat{\vartheta} - \vartheta) + o_p(n^{-1/2})
$$

$$
= O_p(n^{-1/2}),
$$

$$
(3.4) \quad \frac{1}{n} \sum_{i=1}^{n} Z_i \hat{\varepsilon}_i^2 = \frac{1}{n} \sum_{i=1}^{n} Z_i \varepsilon_i^2 + o_p(1) = E Z \sigma^2 + o_p(1),
$$

where $\hat{\vartheta}$ is a $\sqrt{n}$ consistent estimator of $\vartheta$ and $\sigma^2 > 0$ the error variance.

LEMMA 3.1. *Suppose that Assumption* 1 *is satisfied and let* $\hat{\vartheta}$ *be a* $\sqrt{n}$ *consistent estimator of* $\vartheta$. *Then* $\max_{1 \leq j \leq n} |\hat{w}_j - 1| = o_p(1)$ *and*

$$
(3.5) \quad \hat{\lambda} = \frac{1}{\sigma^2} \frac{1}{n} \sum_{j=1}^{n} \frac{Z_j}{EZ} \varepsilon_j - \frac{1}{\sigma^2} E \{\dot{r}_\vartheta(X) | Z = 1\}^\top (\hat{\vartheta} - \vartheta) + o_p(n^{-1/2})
$$

$$
= O_p(n^{-1/2}).
$$

PROOF. We first derive the order of $\hat{\lambda}$. Recall that $\hat{w}_j = 1/(1 + \hat{\lambda} Z_j \hat{\varepsilon}_j)$, that $\hat{w}_j + \hat{\lambda} \hat{w}_j Z_j \hat{\varepsilon}_j = 1$ and that $\sum_{j=1}^{n} \hat{w}_j Z_j \hat{\varepsilon}_j = 0$ by construction. Also note that the $Z_j$'s are binary and that therefore $Z_j = Z_j^2$. This allows us to write

$$
\frac{1}{n} \sum_{j=1}^{n} Z_j \hat{\varepsilon}_j = \frac{1}{n} \sum_{j=1}^{n} (\hat{w}_j + \hat{\lambda} \hat{w}_j Z_j \hat{\varepsilon}_j) Z_j \hat{\varepsilon}_j = \hat{\lambda} \frac{1}{n} \sum_{j=1}^{n} \hat{w}_j Z_j \hat{\varepsilon}_j^2
$$

$$
= \hat{\lambda} \frac{1}{n} \sum_{j=1}^{n} \frac{Z_j \hat{\varepsilon}_j^2}{1 + \hat{\lambda} Z_j \hat{\varepsilon}_j}.
$$

Note that $1 + \hat{\lambda} Z_j \hat{\varepsilon}_j > 0$ since the weights are positive. Then

$$
|\hat{\lambda}| \frac{1}{n} \sum_{j=1}^{n} Z_j \hat{\varepsilon}_j^2 = |\hat{\lambda}| \frac{1}{n} \sum_{j=1}^{n} \frac{Z_j \hat{\varepsilon}_j^2}{1 + \hat{\lambda} Z_j \hat{\varepsilon}_j} (1 + \hat{\lambda} Z_j \hat{\varepsilon}_j)
$$

$$
\leq |\hat{\lambda}| \frac{1}{n} \sum_{j=1}^{n} \frac{Z_j \hat{\varepsilon}_j^2}{1 + \hat{\lambda} Z_j \hat{\varepsilon}_j} \left(1 + |\hat{\lambda}| \max_{1 \leq j \leq n} |Z_j \hat{\varepsilon}_j| \right)
$$

$$
= |\hat{\lambda}| \frac{1}{\hat{\lambda}} \frac{1}{n} \sum_{j=1}^{n} Z_j \hat{\varepsilon}_j \left(1 + |\hat{\lambda}| \max_{1 \leq j \leq n} |Z_j \hat{\varepsilon}_j| \right).
$$

The last equality holds due to (3.6). Applying (3.2), (3.3) and (3.4) to the first and last terms of the inequality we obtain $|\hat{\lambda}| \cdot O_p(1) = O_p(n^{-1/2}) + |\hat{\lambda}| o_p(1)$ which implies $\hat{\lambda} = O_p(n^{-1/2})$. This and (3.2) give $\max_{1 \leq j \leq n} |\hat{\lambda} Z_j \hat{\varepsilon}_j|$

$= o_p(1)$ and therefore our first statement,

$$\max_{1 \le j \le n} |\hat{w}_j - 1| = \max_{1 \le j \le n} \left| \frac{-\hat{\lambda} Z_j \hat{\varepsilon}_j}{1 + \hat{\lambda} Z_j \hat{\varepsilon}_j} \right| = o_p(1).$$

We now again make use of (3.6) and write

$$\frac{1}{n} \sum_{j=1}^{n} Z_j \hat{\varepsilon}_j = \hat{\lambda} \left\{ \frac{1}{n} \sum_{j=1}^{n} (\hat{w}_j - 1) Z_j \hat{\varepsilon}_j^2 + \frac{1}{n} \sum_{j=1}^{n} Z_j \hat{\varepsilon}_j^2 \right\} = \hat{\lambda} \frac{1}{n} \sum_{j=1}^{n} Z_j \hat{\varepsilon}_j^2 + o_p(n^{-1/2}).$$

For the last statement we utilized (3.4), $\max_{1 \le j \le n} |\hat{w}_j - 1| = o_p(1)$ and $\hat{\lambda} = O_p(n^{-1/2})$. This and (3.4) give

$$\hat{\lambda} = \frac{\sum_{j=1}^{n} Z_j \hat{\varepsilon}_j}{\sum_{j=1}^{n} Z_j \hat{\varepsilon}_j^2} + o_p(n^{-1/2})$$

$$= \frac{1}{EZ} \frac{1}{\sigma^2} \frac{1}{n} \sum_{j=1}^{n} Z_j \hat{\varepsilon}_j + o_p(n^{-1/2}).$$

Inserting approximation (3.3) for $n^{-1} \sum_{j=1}^{n} Z_j \hat{\varepsilon}_j$ finally yields the desired approximation of $\hat{\lambda}$. $\square$

LEMMA 3.2.  *Suppose that Assumptions* 1 *and* 2 *are satisfied and let* $\hat{\vartheta}$ *be a* $\sqrt{n}$ *consistent estimator of* $\vartheta$. *Then, with* $\bar{h}(\varepsilon) = E\{h(X, Y)|\varepsilon\}$,

$$\frac{\hat{\lambda}}{n^2} \sum_{i=1}^{n} \sum_{j=1}^{n} \hat{w}_j \hat{\varepsilon}_j \frac{Z_j}{\bar{Z}} h\{X_i, r_{\hat{\vartheta}}(X_i) + \hat{\varepsilon}_j\}$$

$$= \frac{1}{\sigma^2} \frac{1}{n} \sum_{i=1}^{n} \frac{Z_i}{EZ} \varepsilon_i E\{\varepsilon \bar{h}(\varepsilon)\} - \frac{1}{\sigma^2} E\{\varepsilon \bar{h}(\varepsilon)\} E\{\dot{r}_\vartheta(X)|Z = 1\}^\top (\hat{\vartheta} - \vartheta)$$

$$+ o_p(n^{-1/2}).$$

PROOF.   Since $\hat{\lambda} = O_p(n^{-1/2})$ and $\max_{1 \le j \le n} |\hat{w}_j - 1| = o_p(1)$ by the previous lemma, and since $\max_{1 \le i \le n} |Z_i \hat{\varepsilon}_i| = o_p(n^{1/2})$ by (3.2), it is clear that the terms of the sum where $j = i$, that is, $h\{X_i, r_{\hat{\vartheta}}(X_i) + \hat{\varepsilon}_i\} = h(X_i, Y_i)$, can be ignored. It therefore suffices to prove the statement for

$$\frac{\hat{\lambda}}{n^2} \sum_{i} \sum_{j \ne i} \hat{w}_j \hat{\varepsilon}_j \frac{Z_j}{\bar{Z}} h\{X_i, r_{\hat{\vartheta}}(X_i) + \hat{\varepsilon}_j\}$$

$$= \hat{\lambda} \frac{EZ}{\bar{Z}} \psi_w(\hat{\vartheta})$$

with

$$\psi_w(\hat{\vartheta}) = \frac{1}{n^2} \sum_i \sum_{j \neq i} \hat{w}_j \hat{\varepsilon}_j \frac{Z_j}{EZ} h\{X_i, r_{\hat{\vartheta}}(X_i) + \hat{\varepsilon}_j\}$$

$$= \psi(\hat{\vartheta}) + \frac{1}{n^2} \sum_i \sum_{j \neq i} (\hat{w}_j - 1) \hat{\varepsilon}_j \frac{Z_j}{EZ} h\{X_i, r_{\hat{\vartheta}}(X_i) + \hat{\varepsilon}_j\},$$

where $\psi$ is $\psi_w$ with $\hat{w}_j = 1$. The second part involving the difference $\hat{w}_j - 1$ is $o_p(n^{-1/2})$, which can be seen as follows: using $\hat{\lambda} = O_p(n^{-1/2})$ and $\max_{1 \leq j \leq n} |\hat{w}_j - 1| = o_p(1)$ we obtain

$$\left| \hat{\lambda} \frac{EZ}{\bar{Z}} \frac{1}{n^2} \sum_i \sum_{j \neq i} (\hat{w}_j - 1) \hat{\varepsilon}_j \frac{Z_j}{EZ} h\{X_i, r_{\hat{\vartheta}}(X_i) + \hat{\varepsilon}_j\} \right|$$

$$\leq |\hat{\lambda}| \frac{1}{\bar{Z}} \max_{1 \leq j \leq n} |\hat{w}_j - 1| \frac{1}{n^2} \sum_i \sum_{j \neq i} |\hat{\varepsilon}_j h\{X_i, r_{\hat{\vartheta}}(X_i) + \hat{\varepsilon}_j\}|$$

$$= o_p(n^{-1/2}) \cdot \frac{1}{n^2} \sum_i \sum_{j \neq i} |\hat{\varepsilon}_j h\{X_i, r_{\hat{\vartheta}}(X_i) + \hat{\varepsilon}_j\}|.$$

This gives the claimed rate $o_p(n^{-1/2})$ since the sum is bounded in probability, which follows from the $\sqrt{n}$ consistency of $\hat{\vartheta}$ and Assumptions 1 and 2 on the terms of the product $(Y_2 - r_\tau(X_2))h\{X_1, r_\tau(X_1) + Y_2 - r_\tau(X_2)\}$.

It remains to consider $\hat{\lambda} EZ/\bar{Z} \psi(\hat{\vartheta})$. Using $\hat{\lambda} = O_p(n^{-1/2})$ we can replace $\psi(\hat{\vartheta})$ by $\psi(\vartheta)$ since $\psi(\hat{\vartheta}) - \psi(\vartheta) = o_p(1)$, which again follows from Assumptions 1 and 2 and the consistency of $\hat{\vartheta}$. Further, by the law of large numbers, $EZ/\bar{Z} = 1 + o_p(1)$ and $\psi(\vartheta) - E\psi(\vartheta) = o_p(1)$. These arguments yield

$$\hat{\lambda} \frac{EZ}{\bar{Z}} \psi(\hat{\vartheta}) = \hat{\lambda} E\psi(\vartheta) + o_p(n^{-1/2}).$$

The expected value of $\psi(\vartheta)$ is

$$\frac{n-1}{n} E\left[ \varepsilon_2 \frac{Z_2}{EZ} h\{X_1, r_\vartheta(X_1) + \varepsilon_2\} \right] = \frac{n-1}{n} E\{\varepsilon h(X, Y)\} = \frac{n-1}{n} E\{\varepsilon \bar{h}(\varepsilon)\}.$$

Summing up,

$$\hat{\lambda} \frac{EZ}{\bar{Z}} \psi_w(\hat{\vartheta}) = \hat{\lambda} E\psi(\vartheta) + o_p(n^{-1/2}) = \hat{\lambda} E\{\varepsilon \bar{h}(\varepsilon)\} + o_p(n^{-1/2}).$$

Inserting expansion (3.5) for $\hat{\lambda}$ into the above completes the proof. □

Combining the previous lemma and the approximation of the weighted estimator from Section 2 gives an expansion for the weighted estimator.

THEOREM 3.3. *Suppose that Assumption 1 and 2 are satisfied and that $\hat{\vartheta}$ is a $\sqrt{n}$ consistent estimator of $\vartheta$. Let $\bar{h}(\varepsilon) = E\{h(X,Y)|\varepsilon\}$. Then*

$$\frac{1}{n}\sum_{i=1}^{n}\hat{\chi}_w(X_i,\hat{\vartheta}) = \frac{1}{n}\sum_{i=1}^{n}\left(\chi(X_i,\vartheta) + \frac{Z_i}{EZ}\left[\bar{h}(\varepsilon_i) - E\bar{h}(\varepsilon_i) - \frac{E\{\varepsilon\bar{h}(\varepsilon)\}}{\sigma^2}\varepsilon_i\right]\right)$$
$$+ D_w^{\top}(\hat{\vartheta} - \vartheta) + o_p(n^{-1/2}),$$

*where* $D_w = E(h(X,Y)[\dot{r}_\vartheta(X) - E\{\dot{r}_\vartheta(X)|Z = 1\}]\ell(\varepsilon)) + \sigma^{-2}E\{\varepsilon\bar{h}(\varepsilon)\}\times E\{\dot{r}_\vartheta(X)|Z = 1\}$.

PROOF. Consider the two terms of representation (3.1) and replace them by their approximations given in Corollary 2.2 and Lemma 3.2. This yields

$$\frac{1}{n}\sum_{i=1}^{n}\hat{\chi}_w(X_i,\hat{\vartheta})$$
$$= \frac{1}{n}\sum_{i=1}^{n}\left(\chi(X_i,\vartheta) + \frac{Z_i}{EZ}\left[\bar{h}(\varepsilon_i) - E\bar{h}(\varepsilon) - \frac{E\{\varepsilon\bar{h}(\varepsilon)\}}{\sigma^2}\varepsilon_i\right]\right)$$
$$+ \left[D^{\top} + \frac{1}{\sigma^2}E\{\varepsilon\bar{h}(\varepsilon)\}E\{\dot{r}_\vartheta(X)|Z = 1\}^{\top}\right](\hat{\vartheta} - \vartheta) + o_p(n^{-1/2})$$

with $D + \sigma^{-2}E\{\varepsilon\bar{h}(\varepsilon)\}E\{\dot{r}_\vartheta(X)|Z = 1\} = D_w$, by definition of $D$ (see Corollary 2.2). Inserting this into the above gives the desired representation. □

**4. Efficiency.** We are interested in efficient estimation of $Eh(X,Y)$ based on observations $(X, ZY, Z)$. Our estimator requires an efficient estimator of $\vartheta$. In this section we determine the influence function of an efficient estimator of $Eh(X,Y)$. In the next section, where the influence function of an efficient estimator $\hat{\vartheta}$ of $\vartheta$ is determined, this allows us to show that the fully imputed estimator with an efficient $\hat{\vartheta}$ plugged in is efficient. Throughout we will suppose that the assumptions made earlier are satisfied.

We first calculate the efficient influence function for estimating an arbitrary functional $\kappa$ of the joint distribution $P(dx, dy, dz)$. The joint distribution depends on the marginal distribution $G(dx)$ of $X$, the conditional probability $\pi(x)$ of $Z = 1$ given $X = x$, and the conditional distribution $Q(x, dy)$ of $Y$ given $X = x$,

$$P(dx, dy, dz) = G(dx)B_{\pi(x)}(dz)\{zQ(x, dy) + (1 - z)\delta_0(dy)\}.$$

Here $B_p = p\delta_1 + (1-p)\delta_0$ denotes the Bernoulli distribution with parameter $p$ and $\delta_t$ the Dirac measure at $t$. In a first step we consider a nonparametric model for $P$, that is, we allow for arbitrary models for $G, Q$ and $\pi$. For this general setting a characterization of efficient estimators of $\kappa(G, Q, \pi)$ is in Müller, Schick and Wefelmeyer (2006), Section 2. In the following we summarize their key arguments and apply them to the special case of nonlinear regression (which is not considered in that article). We then calculate the efficient influence functions for estimating $Eh(X, Y)$ in the nonlinear regression model and, in the next section, for estimating $\vartheta$.

For the characterization of efficient estimators it is essential to first introduce the notion of tangent spaces. The tangent space of a model is the set of possible perturbations of $P$ within the model. An estimator of a certain functional is, roughly speaking, efficient if its influence function equals the so-called canonical gradient of the functional, which is an element of the tangent space. Hence, in order to characterize the efficient influence function, we first need to determine the tangent space.

Consider (Hellinger differentiable) perturbations of $G$, $Q$ and $\pi$,

$$G_{nu}(dx) \doteq G(dx)\{1 + n^{-1/2}u(x)\},$$

$$Q_{nv}(x, dy) \doteq Q(x, dy)\{1 + n^{-1/2}v(x, y)\},$$

$$B_{\pi_{nw}(x)}(dz) \doteq B_{\pi(x)}(dz)[1 + n^{-1/2}\{z - \pi(x)\}w(x)].$$

To guarantee that the perturbed distributions are probability distributions requires that the (Hellinger) derivative $u$ belongs to

$$L_{2,0}(G) = \left\{u \in L_2(G) \colon \int u \, dG = 0\right\},$$

that $v$ belongs to

$$V_0 = \left\{v \in L_2(M) \colon \int v(x, y) \, Q(x, dy) = 0\right\}$$

with $M(dx, dy) = Q(x, dy)G(dx)$, and that $w$ belongs to $L_2(G_\pi)$, where $G_\pi(dx) = \pi(x)\{1 - \pi(x)\}G(dx)$. The perturbed joint distribution $P_{nuvw}$ then has derivative $t_{uvw}(x, zy, z) = u(x) + zv(x, y) + \{z - \pi(x)\}w(x)$. Note that models for $G, Q$ and $\pi$ will result in further restrictions on the perturbations which must satisfy the model assumptions. Then $u, v$ and $\pi$ must be restricted to subspaces $U$ of $L_{2,0}(G)$, $V$ of $V_0$ and $W$ of $L_2(G_\pi)$.

In this article we make no model assumptions on $G$ and $\pi$ and thus have $U = L_{2,0}(G)$ and $W = L_2(G_\pi)$. Since we are considering nonlinear regression we do, however, have a model for the conditional distribution, namely $Q(x, dy) = f\{y - r_\vartheta(x)\} \, dy$ with $f$ denoting the (mean zero) density of the error distribution. Perturbations $v$ of $Q$ must therefore satisfy $\int v(x, y)f\{y - r_\vartheta(x)\} \, dy = 0$. In order to derive an explicit form of $V$, we introduce perturbations $s$ and $t$ of the two parameters $f$ and $\vartheta$. Write $F$ for the distribution function of $f$ and remember that we assume that $f$ has finite Fisher information for location, $E\ell^2(\varepsilon) < \infty$, where $\ell = -f'/f$ is the score function. The perturbed distribution $Q$ now depends on $s$ and $t$,

$$Q_{nv}(x, dy) = Q_{nst}(x, dy) = f_{ns}\{y - r_{\vartheta_{nt}}(x)\} \, dy$$

with $\vartheta_{nt} = \vartheta + n^{-1/2}t$, $t \in \mathbb{R}^p$, $f_{ns}(y) = f(y)\{1 + n^{-1/2}s(y)\}$ and $s \in S$, where

$$S = \left\{ s \in L_2(F) \colon \int s(y)f(y)\,dy = 0, \int y s(y)f(y)\,dy = 0 \right\}.$$

Note that the space $S$ is determined by two constraints: the perturbed error density $f_{ns}$ must integrate to 1, $\int f_{ns}(y)\,dy = 1$, and must be centered at zero, $\int y f_{ns}(y)\,dy = 0$. As in Schick ([1993](#)), Section 3, we have

$$
\begin{aligned}
f_{ns}&\{y - r_{\vartheta_{nt}}(x)\} \\
&= f\{y - r_{\vartheta_{nt}}(x)\}[1 + n^{-1/2}s\{y - r_{\vartheta_{nt}}(x)\}] \\
&\doteq [f\{y - r_\vartheta(x)\} - n^{-1/2}f'\{y - r_\vartheta(x)\}\dot{r}_\vartheta(x)^\top t][1 + n^{-1/2}s\{y - r_\vartheta(x)\}] \\
&\doteq f\{y - r_\vartheta(x)\}\left(1 + n^{-1/2}\left[s\{y - r_\vartheta(x)\} - \frac{f'\{y - r_\vartheta(x)\}}{f\{y - r_\vartheta(x)\}}\dot{r}_\vartheta(x)^\top t\right]\right) \\
&= f\{y - r_\vartheta(x)\}(1 + n^{-1/2}[s\{y - r_\vartheta(x)\} + \ell\{y - r_\vartheta(x)\}\dot{r}_\vartheta(x)^\top t]).
\end{aligned}
$$

Therefore

$$
\begin{aligned}
Q_{nst}(x, dy) &\doteq f\{y - r_\vartheta(x)\}\,dy \\
&\qquad \times (1 + n^{-1/2}[s\{y - r_\vartheta(x)\} + \ell\{y - r_\vartheta(x)\}\dot{r}_\vartheta(x)^\top t])
\end{aligned}
$$

and the subspace $V$ of $V_0$ is

$$(4.1) \quad V = \{v(x, y) = s\{y - r_\vartheta(x)\} + \ell\{y - r_\vartheta(x)\}\dot{r}_\vartheta(x)^\top t \colon s \in S, t \in \mathbb{R}^p\}.$$

We now briefly review some definitions. We will do this for arbitrary subspaces $U, V$ and $W$ of $L_{2,0}(G)$, $V_0$ and $L_2(G_\pi)$, and then return to our specific situation.

Let $T$ denote the *tangent space* consisting of all derivatives $t_{uvw}$. A functional $\kappa$ of $G$, $Q$ and $\pi$ is called *differentiable* with *gradient* $g \in L_2(P)$ if, for all $u \in U$, $v \in V$ and $w \in W$,

$$
\begin{aligned}
(4.2) \qquad n^{1/2}&\{\kappa(G_{nu}, Q_{nv}, \pi_{nw}) - \kappa(G, Q, \pi)\} \\
&\to E\{g(X, ZY, Z)t_{uvw}(X, ZY, Z)\}.
\end{aligned}
$$

The (unique) *canonical* gradient $g_* = g_*(X, ZY, Z)$ is the projection of $g(X, ZY, Z)$ onto the tangent space $T$. It is easy to check that $T$ can be written as an orthogonal sum of three subspaces,

$$T = \{u(X) \colon u \in U\} \oplus \{Zv(X, Y) \colon v \in V\} \oplus \{\{Z - \pi(X)\}w(X) \colon w \in W\}.$$

The random variable $g_*(X, ZY, Z)$ is therefore the sum $u_*(X) + Zv_*(X, Y) + \{Z - \pi(X)\}w_*(X)$, where $u_*(X)$, $Zv_*(X, Y)$ and $\{Z - \pi(X)\}w_*(X)$ are the projections of $g(X, ZY, Z)$ onto these subspaces.

An estimator $\hat{\kappa}$ for $\kappa$ is *regular* with *limit* $L$ if $L$ is a random variable such that for all $u \in U$, $v \in V$ and $w \in W$,

$$n^{1/2}\{\hat{\kappa} - \kappa(G_{nu}, Q_{nv}, \pi_{nw})\} \Rightarrow L \qquad \text{under } P_{nuvw}.$$

The Hájek–Le Cam convolution theorem says that $L$ is distributed as the sum of a normal random variable $N$, with mean zero and variance $Eg_*^2$, and some independent random variable. This justifies calling an estimator $\hat{\kappa}$ *efficient* if it is regular with limit $L = N$. As a consequence, a regular estimator is efficient if and only if it is asymptotically linear with influence function $g_*$, that is,

$$n^{1/2}\{\hat{\kappa} - \kappa(G, Q, \pi)\} = n^{-1/2} \sum_{i=1}^{n} g_*(X_i, Z_i Y_i, Z_i) + o_p(1).$$

A reference for the convolution theorem and the characterization is Bickel et al. (1998).

Let us now specify the canonical gradient for the functional $Eh(X, Y)$. The canonical gradient is, in particular, a gradient and thus specified by (4.2). Moreover, it is characterized by $g_*(X, ZY, Z) = u_*(X) + Zv_*(X, Y) + \{Z - \pi(X)\}w_*(X)$ with the terms of the sum being projections as stated above. The canonical gradient for arbitrary $\kappa$ is therefore determined by

$$(4.3) \qquad \begin{aligned} &E\{u_*(X)u(X)\} + E\{Zv_*(X, Y)v(X, Y)\} \\ &\quad + E[\{Z - \pi(X)\}^2 w_*(X)w(X)] \\ &\quad = \lim_{n \to \infty} n^{1/2}\{\kappa(G_{nu}, Q_{nv}, \pi_{nw}) - \kappa(G, Q, \pi)\}. \end{aligned}$$

In the nonlinear regression model we have, as defined earlier, $U = L_{2,0}(G)$, $W = L_2(G_\pi)$, $Q_{nv} = Q_{nst}$ with $v \in V$, that is, $v(X, Y) = s(\varepsilon) + \ell(\varepsilon)\dot{r}_\vartheta(X)^\top t$ [see (4.1)]. Since $Eh(X, Y)$ does not depend on $\pi$ we have $Eh(X, Y) = \kappa(G, Q, \pi) = \kappa(G, Q)$ and

$$\begin{aligned} Eh(X, Y) &= \int h \, dM = \iint h(x, y) \, Q(x, dy) G(dx) \\ &= \iint h(x, y) f\{y - r_\vartheta(x)\} \, dy \, G(dx). \end{aligned}$$

Let $M_{nuv}(dx, dy) = Q_{nv}(x, dy) G_{nu}(dx)$ with $Q_{nv} = Q_{nst} = f_{ns}\{y - r_{\vartheta_{nt}}(x)\} \, dy$ and perturbations $G_{nu}$, $f_{ns}$ and $\vartheta_{nt}$ as defined earlier. Using the previous approximations we see that the right-hand side of (4.3) is

$$\lim_{n \to \infty} n^{1/2} \left( \int h \, dM_{nuv} - \int h \, dM \right) = E[h(X, Y)\{u(X) + v(X, Y)\}]$$

with $v(X,Y) = s(\varepsilon) + \ell(\varepsilon)\dot{r}_\vartheta(X)^\top t$. The canonical gradient $g_*$ of $Eh(X,Y)$ is therefore determined by

$$
\begin{aligned}
(4.4) \quad & E\{u_*(X)u(X)\} + E\{Zv_*(X,Y)v(X,Y)\} \\
& \quad + E[\{Z - \pi(X)\}^2 w_*(X)w(X)] = E[h(X,Y)\{u(X) + v(X,Y)\}]
\end{aligned}
$$

for all $u \in U$, $v \in V$ and $w \in W$ with $v$ of the above form.

In order to specify $g_*$ we set $u = 0$ and $v = 0$ in (4.4) and see that $w_*$ must be zero. Setting $v = 0$, we see that $u_*(X)$ is the projection of $h(X,Y)$ onto $U = L_{2,0}(G)$, that is, $u_*(X) = \chi(X,\vartheta) - E\{\chi(X,\vartheta)\}$ with $\chi(X,\vartheta) = E\{h(X,Y)|X\}$. Hence we have

$$
(4.5) \qquad g_*(X, ZY, Z) = \chi(X,\vartheta) - E\{\chi(X,\vartheta)\} + Zv_*(X,Y)
$$

and are left to determine $v_*$. Taking $u = 0$ in (4.4), we see that the projection of $Zv_*(X,Y)$ onto $\tilde{V} = \{v(X,Y): v \in V\}$ must equal the projection of $h(X,Y)$ onto $\tilde{V}$, that is, onto

$$
\tilde{V} = \{s(\varepsilon) + \ell(\varepsilon)\dot{r}_\vartheta(X)^\top t, s \in S, t \in \mathbb{R}^p\}.
$$

There are two possible ways to obtain $v_*$. One method would be to make an educated guess: in Theorem 3.3 we derived an approximation of an estimator of $Eh(X,Y)$ which we expect to be efficient since it uses all information about the model. The approximation still involves $\hat{\vartheta} - \vartheta$ but, combined with the efficient influence function for estimating $\vartheta$ (which is relatively easy to derive; see Section 5), it will suggest a candidate for $v_*$. Whether this candidate is the correct $v_*$ can be checked with characterization (4.4), that is, with

$$
\begin{aligned}
& E[Zv_*(X,Y)\{s(\varepsilon) + \ell(\varepsilon)\dot{r}_\vartheta(X)^\top t\}] = E[h(X,Y)\{s(\varepsilon) + \ell(\varepsilon)\dot{r}_\vartheta(X)^\top t\}]. \\
(4.6) &
\end{aligned}
$$

The other method uses the structure of the tangent space. The canonical gradient $v_*$ is characterized in terms of projections onto $\tilde{V}$. Its derivation as a projection onto $\tilde{V}$ is simplified by decomposing $\tilde{V}$. Let $\ell_s$ denote the projection of $\ell$ onto $S$,

$$
\ell_s(\varepsilon) = \ell(\varepsilon) - \sigma^{-2}\varepsilon,
$$

and note that $\ell_s = 0$ is possible, namely when the error density $f$ is normal. We now introduce the notation

$$
\zeta = [\dot{r}_\vartheta(X) - E\{\dot{r}_\vartheta(X)|Z = 1\}]\ell(\varepsilon) + E\{\dot{r}_\vartheta(X)|Z = 1\}\frac{\varepsilon}{\sigma^2}
$$

and, for $s \in S$ and $t \in \mathbb{R}^p$, write

$$
\begin{aligned}
& s(\varepsilon) + \dot{r}_\vartheta(X)^\top t\ell(\varepsilon) \\
& \quad = s(\varepsilon) + t^\top[\dot{r}_\vartheta(X) - E\{\dot{r}_\vartheta(X)|Z = 1\}]\ell(\varepsilon)
\end{aligned}
$$

$$+ t^\top E\{\dot{r}_\vartheta(X)|Z=1\}\Big\{\ell(\varepsilon) - \frac{\varepsilon}{\sigma^2}\Big\} + t^\top E\{\dot{r}_\vartheta(X)|Z=1\}\frac{\varepsilon}{\sigma^2}$$

$$= t^\top \zeta + s(\varepsilon) + t^\top E\{\dot{r}_\vartheta(X)|Z=1\}\ell_s(\varepsilon)$$

with $s(\varepsilon) + t^\top E\{\dot{r}_\vartheta(X)|Z=1\}\ell_s(\varepsilon) \in S$. Any element of $\tilde{V}$ can therefore be written $t^\top \zeta + s(\varepsilon)$ for some $t \in \mathbb{R}^p$ and $s \in S$. Since the canonical gradient $v_*$ is in $\tilde{V}$ by definition, it must be of the form

$$v_*(X,Y) = s^*(\varepsilon) + t^{*\top}\zeta$$

with $s^* \in S$ and $t^* \in \mathbb{R}^p$ to be determined such that (4.6) holds, that is, after our above considerations,

$$E[Z\{s^*(\varepsilon) + t^{*\top}\zeta\}\{s(\varepsilon) + t^\top\zeta\}] = E[h(X,Y)\{s(\varepsilon) + t^\top\zeta\}]$$

for all $t \in \mathbb{R}^p$ and $s \in S$.

We first consider $t = 0$ and secondly $s = 0$ and, in both cases, use the fact that $Z\zeta$ is orthogonal to $S$. Then the above characterization of $s^*$ and $t^*$ reduces to two equations, namely

$$(4.7) \qquad E\{Zs^*(\varepsilon)s(\varepsilon)\} = E\{h(X,Y)s(\varepsilon)\} \qquad \text{for all } s \in S,$$

$$(4.8) \qquad E\{Zt^{*\top}\zeta t^\top\zeta\} = E\{h(X,Y)t^\top\zeta\} \qquad \text{for all } t \in \mathbb{R}^p.$$

Consider (4.7) and again use the notation $\bar{h}(\varepsilon)$ for the conditional expectation $E\{h(X,Y)|\varepsilon\}$. Then (4.7) can be written as $E\{Zs^*(\varepsilon)s(\varepsilon)\} = E\{\bar{h}(\varepsilon)s(\varepsilon)\}$, that is, $\bar{h}(\varepsilon)/EZ$ is an obvious candidate for $s^*$. However, it is not (yet) in $S$: the desired $s^*$ is obtained as its centered version with a correction term chosen such that $s^* \in S$,

$$s^*(\varepsilon) = \frac{1}{EZ}\Big[\bar{h}(\varepsilon) - E\bar{h}(\varepsilon) - \frac{E\{\varepsilon\bar{h}(\varepsilon)\}}{\sigma^2}\varepsilon\Big].$$

The vector $t^*$ is obtained by solving (4.8), $t^{*\top}E(Z\zeta\zeta^\top)t = E\{h(X,Y)\zeta^\top\}t$ for all $t \in \mathbb{R}^p$. Now use the definition of $\zeta$ from above and the definition of the vector $D_w$ from the end of the previous section, $D_w = E(h(X,Y)[\dot{r}_\vartheta(X) - E\{\dot{r}_\vartheta(X)|Z=1\}]\ell(\varepsilon)) + \sigma^{-2}E\{\varepsilon\bar{h}(\varepsilon)\}E\{\dot{r}_\vartheta(X)|Z=1\}$, and assume that $E(Z\zeta\zeta^\top)$ is invertible to obtain

$$t^{*\top} = E\{h(X,Y)\zeta^\top\}E(Z\zeta\zeta^\top)^{-1}$$

$$= E\Big\{h(X,Y)\Big([\dot{r}_\vartheta(X) - E\{\dot{r}_\vartheta(X)|Z=1\}]^\top\ell(\varepsilon)$$

$$\qquad\qquad + E\{\dot{r}_\vartheta(X)|Z=1\}^\top\frac{\varepsilon}{\sigma^2}\Big)\Big\}E(Z\zeta\zeta^\top)^{-1}$$

$$= D_w^\top E(Z\zeta\zeta^\top)^{-1}.$$

This completes the derivation of $v_*(X, Y) = s^*(\varepsilon) + t^{*\top}\zeta$:

$$(4.9) \quad v_*(X, Y) = \frac{1}{EZ}\left[\bar{h}(\varepsilon) - E\bar{h}(\varepsilon) - \frac{E\{\varepsilon\bar{h}(\varepsilon)\}}{\sigma^2}\varepsilon\right] + D_w^\top E(Z\zeta\zeta^\top)^{-1}\zeta.$$

Equations (4.5) and (4.9) together finally yield the canonical gradient $g_*$, which is given in the following lemma. Note that we now have the additional assumption that $E(Z\zeta\zeta^\top)$ is invertible, where $E(Z\zeta\zeta^\top)$ involves the covariance matrix of $Z\dot{r}(X)$ and the Fisher information $E\dot\ell^2(\varepsilon)$.

LEMMA 4.1. *Let* $\bar{h}(\varepsilon) = E\{h(X, Y)|\varepsilon\}$, $\zeta = [\dot{r}_\vartheta(X) - E\{\dot{r}_\vartheta(X)|Z=1\}]\ell(\varepsilon) + \sigma^{-2}E\{\dot{r}_\vartheta(X)|Z=1\}\varepsilon$ *and* $D_w = E(h(X, Y)[\dot{r}_\vartheta(X) - E\{\dot{r}_\vartheta(X)|Z=1\}]\ell(\varepsilon)) + \sigma^{-2}E\{\varepsilon\bar{h}(\varepsilon)\}E\{\dot{r}_\vartheta(X)|Z=1\} = E\{h(X, Y)\zeta\}$. *Suppose additionally to the model assumptions from Section 2 that* $E(Z\zeta\zeta^\top)$ *is invertible. Then the canonical gradient of the functional* $Eh(X, Y)$ *is*

$$(4.10) \quad \begin{aligned} &\chi(X, \vartheta) - E\{\chi(X, \vartheta)\} \\ &\quad + \frac{Z}{EZ}\left[\bar{h}(\varepsilon) - E\bar{h}(\varepsilon) - \frac{E\{\varepsilon\bar{h}(\varepsilon)\}}{\sigma^2}\varepsilon\right] + D_w^\top E\{Z\zeta\zeta^\top\}^{-1}Z\zeta. \end{aligned}$$

## 5. Estimation of the parameter and main result.

In this section we show that the weighted estimator for $Eh(X, Y)$ with an efficient estimator $\hat\vartheta$ for $\vartheta$ plugged in is asymptotically linear with influence function equal to the canonical gradient, that is, it is efficient. Let us compare the expansion of the weighted estimator from Theorem 3.3 and the efficient influence function which is given by the canonical gradient (4.10) in Lemma 4.1. The approximation of $n^{-1/2}\sum_{i=1}^n[\hat\chi_w(X_i, \hat\vartheta) - E\{\chi(X, \vartheta)\}]$ which we derived in Section 3 is

$$\begin{aligned} n^{-1/2}\sum_{i=1}^n\Big(&\chi(X_i, \vartheta) - E\{\chi(X, \vartheta)\} \\ &+ \frac{Z_i}{EZ}\left[\bar{h}(\varepsilon_i) - E\bar{h}(\varepsilon) - \frac{E\{\varepsilon\bar{h}(\varepsilon)\}}{\sigma^2}\varepsilon_i\right]\Big) + D_w^\top n^{1/2}(\hat\vartheta - \vartheta), \end{aligned}$$

where $D_w = E(h(X, Y)[\dot{r}_\vartheta(X) - E\{\dot{r}_\vartheta(X)|Z=1\}]\ell(\varepsilon)) + \sigma^{-2}E\{\varepsilon\bar{h}(\varepsilon)\} \times E\{\dot{r}_\vartheta(X)|Z=1\}$. The efficient influence function determined by the canonical gradient is

$$\begin{aligned} &\chi(X, \vartheta) - E\{\chi(X, \vartheta)\} \\ &\quad + \frac{Z}{EZ}\left[\bar{h}(\varepsilon) - E\bar{h}(\varepsilon) - \frac{E\{\varepsilon\bar{h}(\varepsilon)\}}{\sigma^2}\varepsilon\right] + D_w^\top E\{Z\zeta\zeta^\top\}^{-1}Z\zeta \end{aligned}$$

with $\zeta = [\dot{r}_\vartheta(X) - E\{\dot{r}_\vartheta(X)|Z=1\}]\ell(\varepsilon) + \sigma^{-2}E\{\dot{r}_\vartheta(X)|Z=1\}\varepsilon$. Using an estimator $\hat\vartheta$ with influence function $E(Z\zeta\zeta^\top)^{-1}Z\zeta$ would therefore yield an

efficient estimator for $Eh(X, Y)$. In fact, it is easy to check (this will be done in the following lemma) that this influence function is the canonical gradient of the functional $\kappa(G, Q, \pi) = \vartheta$. This means that our estimator of $Eh(X, Y)$ requires an efficient estimator $\hat{\vartheta}$ for $\vartheta$ to be plugged in in order to be efficient.

**Lemma 5.1.** *Let $\zeta = [\dot{r}_\vartheta(X) - E\{\dot{r}_\vartheta(X)|Z = 1\}]\ell(\varepsilon) + \sigma^{-2}E\{\dot{r}_\vartheta(X)|Z = 1\}\varepsilon$ and suppose that $E(Z\zeta\zeta^\top)$ is invertible. An asymptotically linear estimator $\hat{\vartheta}$ for $\vartheta$ with influence function $E(Z\zeta\zeta^\top)^{-1}Z\zeta$, that is,*

$$n^{1/2}(\hat{\vartheta} - \vartheta)$$
$$= n^{-1/2} \sum_{i=1}^n E(Z\zeta\zeta^\top)^{-1} Z_i \left[ \left\{ \dot{r}_\vartheta(X_i) - E[\dot{r}_\vartheta(X)|Z = 1] \right\} \ell(\varepsilon_i) \right.$$
$$\left. + E\{\dot{r}_\vartheta(X)|Z = 1\} \frac{\varepsilon_i}{\sigma^2} \right] + o_p(1),$$

*is efficient for $\vartheta$.*

**Proof.** We have a semiparametric model for the conditional distribution, namely $Q(x, dy) = f(y - r_\vartheta(x)) \, dy$, and nonparametric models for $G$ and $\pi$. The functional $\vartheta \in \mathbb{R}^p$ is therefore a functional of $Q$, $\kappa(G, Q, \pi) = \kappa(Q) = \vartheta$. By the discussion of the previous section we must show that the influence function of the estimator equals the canonical gradient, which is, for arbitrary functionals $\kappa$, determined by (4.3). For the functional $\vartheta$ the right-hand side of (4.3) is simply $n^{1/2}\{(\vartheta + n^{-1/2}t) - \vartheta\} = t$. From Section 4 we also know that in the nonlinear regression model any $v$ in $\tilde{V}$ is of the form $v(X, Y) = s(\varepsilon) + t^\top \zeta$, where $s \in S$ and $t \in \mathbb{R}$. The canonical gradient $u_*(X) + Zv_*(X, Y) + \{Z - \pi(X)\}w_*(X)$ is therefore characterized by

$$E\{u_*(X)u(X)\} + E[Zv_*(X, Y)\{s(\varepsilon) + \zeta^\top t\}]$$
$$+ E[\{Z - \pi(X)\}^2 w_*(X)w(X)] = t.$$

Taking $s = 0$, $t = 0$ and $w = 0$ we see that $u_* = 0$. Analogously one obtains that $w_*$ must be zero. The canonical gradient thus reduces to $Zv_*(X, Y)$. Again, since $v_* \in \tilde{V}$, we write $Zv_*(X, Y) = Zs^*(\varepsilon) + Z\zeta^\top t^*$ with $s^*$ and $t^*$ to be determined. Taking $t = 0$ we see that $Zv_*$ must be orthogonal to $S$, that is, $s^* = 0$ which yields $Zv_*(X, Y) = Z\zeta^\top t^*$. The above characterization therefore reduces to

$$t = E[Z\zeta^\top t^*\{s(\varepsilon) + \zeta^\top t\}] = t^{*\top} E(Z\zeta\zeta^\top)t \qquad \text{for all } t \in \mathbb{R}.$$

This gives $t^* = E(Z\zeta\zeta^\top)^{-1}$ and the proof is complete: the canonical gradient of the parameter $\vartheta$ is $Zv_*(X, Y) = Zt^{*\top}\zeta = E(Z\zeta\zeta^\top)^{-1}Z\zeta$. $\square$

Note that the asymptotic variance of $\hat{\vartheta}$ is $E(Z\zeta\zeta^{\top})^{-1}$. The assumption that $E(Z\zeta\zeta^{\top})$ must be invertible is therefore a condition on the covariance matrix of an efficient estimator of $\vartheta$ which we require to have full rank. Lemma 5.1 combined with the previous discussion yields our main result, which is given in the following theorem. Note that the asymptotic variance of the fully imputed estimator of $Eh(X,Y)$ is $Eg_*^2$, where $g_*$ is the canonical gradient from (4.10). This variance is also given in the theorem below and is easily verified by taking into account that the three terms of $g_*$ are orthogonal.

THEOREM 5.2. *Assume that Assumptions 1 and 2 hold and that the covariance matrices of $\dot{r}_\vartheta(X)$ and of $Z\dot{r}_\vartheta(X)$ are invertible. Let $\hat{\vartheta}$ be an asymptotically linear estimator of $\vartheta$ with influence function $E(Z\zeta\zeta^{\top})^{-1}Z\zeta$ where $\zeta = [\dot{r}_\vartheta(X) - E\{\dot{r}_\vartheta(X)|Z=1\}]\ell(\varepsilon) + \sigma^{-2}E\{\dot{r}_\vartheta(X)|Z=1\}\varepsilon$. Then the estimator $n^{-1}\sum_{i=1}^n \hat{\chi}_w(X_i, \hat{\vartheta})$ with $\hat{\chi}_w(X_i, \hat{\vartheta}) = \sum_{j=1}^n \hat{w}_j Z_j h\{x, r_{\hat{\vartheta}}(x) + Y_j - r_{\hat{\vartheta}}(X_j)\}/\sum_{j=1}^n Z_j$ has the expansion*

$$\frac{1}{n}\sum_{i=1}^n \Bigg(\chi(X_i, \vartheta) + \frac{Z_i}{EZ}\Big[\bar{h}(\varepsilon_i) - E\bar{h}(\varepsilon_i) - \frac{E\{\varepsilon\bar{h}(\varepsilon)\}}{\sigma^2}\varepsilon_i\Big]$$

$$+ D_w^{\top} E(Z\zeta\zeta^{\top})^{-1} Z_i$$

$$\times [\dot{r}_\vartheta(X_i) - E\{\dot{r}_\vartheta(X)|Z=1\}]\ell(\varepsilon_i) + E\{\dot{r}_\vartheta(X)|Z=1\}\frac{\varepsilon_i}{\sigma^2}\Bigg)$$

$$+ o_p(n^{-1/2}),$$

*where $D_w = E(h(X,Y)[\dot{r}_\vartheta(X) - E\{\dot{r}_\vartheta(X)|Z=1\}]\ell(\varepsilon)) + \sigma^{-2}E\{\varepsilon\bar{h}(\varepsilon)\} \times E\{\dot{r}_\vartheta(X)|Z=1\}$ and $\bar{h}(\varepsilon) = E\{h(X,Y)|\varepsilon\}$. In particular, it is an efficient estimator of $Eh(X,Y)$ and asymptotically normally distributed with asymptotic variance*

$$E\chi^2(X,\vartheta) + \frac{1}{EZ}E\bar{h}^2(\varepsilon) - \Big(1 + \frac{1}{EZ}\Big)E^2h(X,Y) - \frac{E^2\{\varepsilon\bar{h}(\varepsilon)\}}{\sigma^2 EZ}$$

$$+ D_w^{\top} E(Z\zeta\zeta^{\top})^{-1} D_w.$$

In the linear regression model *without* missing responses, efficient estimators for $\vartheta$ have been constructed by Bickel (1982), Koul and Susarla (1983) and Schick (1987, 1993). Schick (1993) considers general regression models with arbitrary sets of identifiability assumptions and discusses the mean zero constraint on the error distribution as an important example. His construction of an efficient estimator requires a preliminary estimate of $\vartheta$ and a direct estimator of the influence function. The influence function for the nonlinear regression model with mean zero errors [see Schick (1993), Section 4.1 and Remark 3.13] is $E(\xi\xi^{\top})^{-1}\xi$ with $\xi = [\dot{r}_\vartheta(X) - E\{\dot{r}_\vartheta(X)\}]\ell(\varepsilon) + E[\dot{r}_\vartheta(X)]\varepsilon/\sigma^2$

and therefore consistent with our findings. A further developed efficient estimator, which requires weaker conditions, is in Forrester et al. 2003. In the model with missing responses an efficient estimator can be constructed analogously, using only the (available) full observations. Note that the only difference in the construction is that the data are incomplete, that is, the presence of indicators $Z_i$. In the following we will briefly sketch this "one-step improvement" construction of the estimator and refer to Forrester et al. 2003 for details.

Let $\bar{\vartheta}$ denote a $\sqrt{n}$ consistent and *discretized* estimator of $\vartheta$, that is, with values on a rectangular grid with side lengths of order $n^{-1/2}$. Write $\mu(\vartheta)$ for $E\{\dot{r}_\vartheta(X)|Z=1\}$, $\varepsilon(\vartheta)$ for the error variables $\varepsilon(\vartheta) = Y - r_\vartheta(X)$ and $\zeta_\vartheta\{X, \varepsilon(\vartheta)\}$ for $\zeta$, that is,

$$\zeta = \zeta_\vartheta\{X, \varepsilon(\vartheta)\} = \{\dot{r}_\vartheta(X) - \mu(\vartheta)\}\ell\{\varepsilon(\vartheta)\} + \mu(\vartheta)\varepsilon(\vartheta)/\sigma^2.$$

In order to estimate the influence function one replaces the unknown quantities by estimators. The estimator of $\vartheta$ is then of the form

$$\bar{\vartheta} + \left[\sum_{j=1}^n Z_j \hat{\zeta}_{\bar{\vartheta}}\{X_j, \varepsilon_j(\bar{\vartheta})\} \hat{\zeta}_{\bar{\vartheta}}\{X_j, \varepsilon_j(\bar{\vartheta})^\top\}\right]^{-1} \sum_{j=1}^n Z_j \hat{\zeta}_{\bar{\vartheta}}\{X_j, \varepsilon_j(\bar{\vartheta})\},$$

where

$$\hat{\zeta}_{\bar{\vartheta}}(X, \varepsilon(\bar{\vartheta})) = [\dot{r}_{\bar{\vartheta}}(X) - \hat{\mu}(\bar{\vartheta})]\hat{\ell}\{\varepsilon(\bar{\vartheta})\} + \hat{\mu}(\bar{\vartheta})\varepsilon(\bar{\vartheta})/\sigma^2(\bar{\vartheta})$$

with

$$\hat{\mu}(\bar{\vartheta}) = \frac{\sum_{j=1}^n Z_j \dot{r}_{\bar{\vartheta}}(X_j)}{\sum_{j=1}^n Z_j}, \qquad \hat{\sigma}^2(\bar{\vartheta}) = \frac{\sum_{j=1}^n Z_j \varepsilon_j(\bar{\vartheta})^2}{\sum_{j=1}^n Z_j}$$

and an estimator $\hat{\ell}$ of the score function. To describe this estimator let $k$ be a kernel that satisfies the assumptions given in Section 8 of Forrester et al., for example, a logistic density. For a bandwidth $a_n \to 0$ we set $k_n(x) = k(x/a_n)/a_n$. The estimator of the score function $\ell$ is a kernel estimator based on the available residuals $\varepsilon(\bar{\vartheta})$,

$$\hat{\ell}_{\bar{\vartheta}}(x) = \frac{-\hat{f}'_n(x)}{b_n + \hat{f}_n(x)}$$

with $\hat{f}_n(x) = n^{-1} \sum_{j=1}^n Z_j k_n\{x - \varepsilon_j(\bar{\vartheta})\}$ and where $b_n$ is a sequence of positive numbers converging to zero. The orders of $a_n \to 0$ and $b_n \to 0$ (which also apply if only a fixed fraction of the $n$ data pairs is observed) are given in Forrester et al. 2003.

There are other simple estimators for $\vartheta$ available which, however, and in contrast to the estimators proposed by Schick (1987, 1993) and Forrester et al. 2003, are not efficient for $\vartheta$ and which, if used for plug-in, would yield inefficient estimators of $Eh(X,Y)$. One could, for example, estimate $\hat{\vartheta}$ by a weighted least squares estimator, that is, by the solution $t = \hat{\vartheta}$ of an estimating equation $\sum_{i=1}^n Z_i w_t(X_i)\{Y_i - r_t(X_i)\} = 0$. Such an estimator would be appropriate in a regression model where independence of errors and covariates cannot be assumed. Then one could even obtain efficiency for suitably chosen weights [see Müller (2007), for nonlinear regression without missing responses]. The estimating equation can be regarded as an empirical version of the equation $E[Zw_t(X)\{Y - r_t(X)\}] = 0$. If a solution $t = \vartheta$ of this equation exists, the solution $\hat{\vartheta}$ of the empirical version will, in general, be consistent for $\vartheta$. If one is not interested in efficiency, the estimator $n^{-1}\sum_{i=1}^n \hat{\chi}_w(X_i, \hat{\vartheta})$ with a least squares estimator $\hat{\vartheta}$ plugged in would yield a consistent estimator for $Eh(X,Y)$ (but not an efficient one since the independence structure is not used). Alternatively, the least squares estimator can be used as a preliminary estimator for the one-step improvement approach sketched above.

**6. Special cases, simulations and inference.** Sometimes the estimator simplifies considerably, especially if we study simple special cases such as estimation of expectations $Eh(X,Y)$ where $h$ has a simple form. The main result from Theorem 5.2 is therefore useful in proving efficiency of existing approaches for specific applications, or in improving them, and for comparisons of competing methods. Theorem 5.2 further provides the limiting distribution of the efficient estimator, which facilitates the construction of confidence intervals. We will address this and aspects of the construction of estimators in the following, and illustrate the results with simulations.

6.1. *Special cases.* We have shown that the fully imputed weighted estimator $n^{-1}\sum_{i=1}^n \hat{\chi}_w(X_i, \hat{\vartheta})$ with

$$\hat{\chi}_w(x, \hat{\vartheta}) = \sum_{j=1}^n \hat{w}_j Z_j h\{x, r_{\hat{\vartheta}}(x) + Y_j - r_{\hat{\vartheta}}(X_j)\} \Big/ \sum_{j=1}^n Z_j$$

is efficient for $Eh(X,Y)$ where $h(X,Y)$ is a known square-integrable function. The literature usually deals with estimation of the mean response, that is, $h(x,y) = y$. Other important examples are estimation of higher moments of the response variable $Y$ and the estimation of the covariance and of mixed moments of $X$ and $Y$. In all these cases $h(x,y)$ is a polynomial in $x$ and $y$ and the estimator often simplifies. This holds for the mean response, and, more generally, when $h$ is of the form $h(x,y) = a(x)y$. Then the estimator reduces to an unweighted empirical estimator, which can be seen as follows. Recall that the weights must be chosen such that $\sum_{j=1}^n \hat{w}_j Z_j \hat{\varepsilon}_j = 0$ and that $\hat{w}_j = 1 - \hat{\lambda}\hat{w}_j Z_j \hat{\varepsilon}_j$ which gives $\sum_{j=1}^n \hat{w}_j Z_j / \sum_{j=1}^n Z_j = 1$. Hence the estimator

for $E\{a(X)Y\}$ is

$$\frac{1}{n}\sum_{i=1}^{n}\hat{\chi}_w(X_i,\hat{\vartheta}) = \frac{1}{n}\sum_{i=1}^{n}\frac{\sum_{j=1}^{n}\hat{w}_j Z_j a(X_i)\{r_{\hat{\vartheta}}(X_i)+\hat{\varepsilon}_j\}}{\sum_{j=1}^{n}Z_j}$$

$$= \frac{1}{n}\sum_{i=1}^{n}\frac{\sum_{j=1}^{n}\hat{w}_j Z_j}{\sum_{j=1}^{n}Z_j}a(X_i)r_{\hat{\vartheta}}(X_i) + \frac{\sum_{j=1}^{n}\hat{w}_j Z_j \hat{\varepsilon}_j}{\sum_{j=1}^{n}Z_j}$$

$$= \frac{1}{n}\sum_{i=1}^{n}a(X_i)r_{\hat{\vartheta}}(X_i).$$

In these cases it is therefore not necessary to determine weights: the above intuitive estimator, with an efficient estimator $\hat{\vartheta}$ for $\vartheta$ plugged in, is efficient for $E\{a(X)Y\}$.

An interesting special case is estimation of the mean response, $a(X) = 1$, when *possibly all* responses are observed, which we mentioned in the Introduction. Regardless of whether there are missing responses or not, $n^{-1}\sum_{i=1}^{n}r_{\hat{\vartheta}}(X_i)$ is efficient for $EY$, provided $\hat{\vartheta}$ is efficient for $\vartheta$. The difference between the two situations is the construction of $\hat{\vartheta}$, which will be based on either complete data pairs or on missing response data. Let us stay with this example and consider, for a comparison, the unweighted estimator (1.1) from the introduction, that is, with all weights equal to one. It involves the term $\sum_{j=1}^{n}Z_j\hat{\varepsilon}_j/\sum_{j=1}^{n}Z_j$ which is nonzero. If all responses are observable, the unweighted estimator further simplifies, namely to

$$\frac{1}{n}\sum_{i=1}^{n}r_{\hat{\vartheta}}(X_i) + \frac{1}{n}\sum_{i=1}^{n}\hat{\varepsilon}_i = \frac{1}{n}\sum_{i=1}^{n}Y_i$$

[whereas the weighted estimator is $n^{-1}\sum_{i=1}^{n}r_{\hat{\vartheta}}(X_i)$]. Its influence function is $Y - EY$ which is clearly not the efficient one: our efficient estimator for $EY$ (with an efficient estimator $\hat{\vartheta}$) has the expansion

$$\frac{1}{n}\sum_{i=1}^{n}r_{\hat{\vartheta}}(X_i) \doteq \frac{1}{n}\sum_{i=1}^{n}r_{\vartheta}(X_i) + (\hat{\vartheta}-\vartheta)E\dot{r}_{\vartheta}(X).$$

We recognize this as the expansion from Theorem 3.3 with $D_w = E\dot{r}_{\vartheta}(X)$. Even without inserting the expansion for $\hat{\vartheta}-\vartheta$ from the previous section, it is clear that this is, in general, not the influence function of $n^{-1}\sum_{i=1}^{n}Y_i$, which shows that it cannot be efficient. Note that $n^{-1}\sum_{i=1}^{n}Y_i$ also coincides with the (inefficient) partially imputed estimator if all responses were observed.

6.2. *Simulations.* For an illustration with computer simulations we consider a linear regression function, $r_{\vartheta}(X) = \vartheta X$ with $\vartheta = 2$, and a nonlinear regression function, $r_{\vartheta}(X) = \cos(\vartheta X)$, also with $\vartheta = 2$. The probabilities $\pi(X) = P(Z = 1|X) = E(Z|X)$ are chosen as values of a logistic distribution function, $\pi(X) = 1/(1 + e^{-X})$, so that on average one half of the simulated responses are missing. We generate covariates $X$ from a uniform distribution on the interval $(-1, 1)$ and error variables $\varepsilon$ from a standard normal distribution. If the errors are in fact normally distributed then $\ell(\varepsilon) = \varepsilon/\sigma^2$ and the efficient one-step improvement estimator for $\vartheta$ from the previous section is asymptotically equivalent to the ordinary least squares estimator. The following considerations can therefore be based on this straightforward estimation approach.

In a first example we consider estimation of the mean response $EY$ and compare the efficient (fully imputed weighted) estimator, which, as seen above, here simplifies to $n^{-1} \sum_{i=1}^{n} r_{\hat{\vartheta}}(X_i)$, with the partially imputed estimator $n^{-1} \sum_{i=1}^{n} \{Z_i Y_i + (1 - Z_i) r_{\hat{\vartheta}}(X_i)\}$. We also study the performance of these estimators if the parameter estimates are replaced by their true values, and if all responses are observed, $\pi(\cdot) = 1$. Further we calculate the first simple estimator from the introduction, $n^{-1} \sum_{i=1}^{n} Z_i Y_i / \hat{\pi}(X_i)$, with, for reasons of simplicity, the estimated probabilities $\hat{\pi}$ replaced by the true ones. The values of the simulated mean squared errors are given in Table 1.

In both the linear and the nonlinear regression models, the fully imputed estimator performs considerably better than the partially imputed estimator. The simple estimator in the last column is clearly outperformed by the imputation approaches. Comparing the columns for the fully imputed estimator with and without parameter estimation (and analogously for the partially imputed estimator), we see that the estimator of the slope $\vartheta$ in linear regression $r_\vartheta(X) = \vartheta X$ is, as a plug-in estimator for estimating $EY$, better than the parameter estimator of the frequency parameter $\vartheta$ in the nonlinear regression model $r_\vartheta(X) = \cos(\vartheta X)$: in the linear regression model the mean squared errors of the approaches based on $\vartheta$ and $\hat{\vartheta}$ are very similar, in contrast to the nonlinear model where the differences are quite large. Let us also compare the (a) and (b) sections in the linear regression and the nonlinear regression example, which refer to the situation where (a) responses are missing at random and (b) all responses are available. For the fully imputed estimator $n^{-1} \sum_{i=1}^{n} r_{\hat{\vartheta}}(X_i)$ we observe the expected improved performance when more (response) data for the estimation of $\vartheta$ are available. The situation is different for the partially imputed estimator. Indeed we expect that, similarly, performance will improve as the proportion of observed responses increases. In this case $\hat{\vartheta}$ improves as an estimator of $\vartheta$ but, at the same time, the partially imputed estimator will discard more and more information about the structure of the regression function. [In the extreme case $\pi(\cdot) = 1$ it equals the empirical estimator $n^{-1} \sum_{i=1}^{n} Y_i$.] Our example demonstrates that both scenarios are possible: for the linear regression model the estimator of $\vartheta$ performs well and the simulated mean

Table 1
*Simulated mean squared errors of estimators of the mean response $EY$*

| $\pi(X)$ | $n$ | $\widehat{\text{FI}}$ | FI | $\widehat{\text{PI}}$ | PI | N |
|---|---|---|---|---|---|---|
| Linear regression: $r_\vartheta(X) = \vartheta X$ $(\vartheta = 2)$ | | | | | | |
| $1/(1+e^{-X})$ | 50 | 0.027520 | 0.026639 | 0.036231 | 0.036368 | 0.104962 |
| | 100 | 0.013502 | 0.013298 | 0.018074 | 0.018364 | 0.052680 |
| | 1000 | 0.001328 | 0.001325 | 0.001794 | 0.001835 | 0.005270 |
| 1 | 50 | 0.026990 | 0.026639 | 0.046322 | 0.046322 | 0.046322 |
| | 100 | 0.013415 | 0.013298 | 0.023479 | 0.023479 | 0.023479 |
| | 1000 | 0.001327 | 0.001325 | 0.002345 | 0.002345 | 0.002345 |
| Nonlinear regression: $r_\vartheta(X) = \cos(\vartheta X)$ $(\vartheta = 2)$ | | | | | | |
| $1/(1+e^{-X})$ | 50 | 0.027858 | 0.003957 | 0.031163 | 0.013272 | 0.053038 |
| | 100 | 0.015462 | 0.002001 | 0.017147 | 0.007020 | 0.028154 |
| | 1000 | 0.001492 | 0.000199 | 0.001671 | 0.000696 | 0.002810 |
| 1 | 50 | 0.016512 | 0.003957 | 0.023369 | 0.023369 | 0.023369 |
| | 100 | 0.008581 | 0.002001 | 0.012043 | 0.012043 | 0.012043 |
| | 1000 | 0.000852 | 0.000199 | 0.001207 | 0.001207 | 0.001207 |

*Notes.* The table entries are the simulated mean squared errors of estimators of $EY = Er_\vartheta(X)$ with partially missing responses, $\pi(X) = 1/(1+e^{-X})$ and completely observed data pairs, $\pi(X) = 1$. In the first two columns we study the efficient fully imputed weighted estimator with the ordinary least squares estimator $\hat\vartheta$ plugged in ($\widehat{\text{FI}}$) and its corresponding version using the true parameter, $\vartheta = 2$ (FI). The next two columns refer to the partially imputed estimator using $\hat\vartheta$ ($\widehat{\text{PI}}$) and the version based on $\vartheta = 2$ (PI). The last column considers the simple estimator $n^{-1}\sum_{i=1}^n Z_i Y_i/\pi(X_i)$ (N), which does *not* use imputation. Note that in the sections with $\pi(X) = 1$ the columns for $\widehat{\text{PI}}$, PI and N are identical: since all the indicators are 1, these estimators coincide with the empirical estimator $n^{-1}\sum_{i=1}^n Y_i$.

squared error of the partially imputed estimator in (a) is smaller than in (b). In the nonlinear regression model the estimator of $\vartheta$ is not as good and the mean squared error in (a) is larger than the mean squared error of the empirical estimator in (b). Note that this observation about the performance of the partially imputed estimator is only of secondary interest since, in any case, the fully imputed estimator has the smaller mean squared error.

The situation is slightly more complicated when $h$ is of the form $h(x,y) = a(x)b(y)$ with a *nonlinear* function $b$, for example, when higher mixed moments of $X$ and $Y$ or just higher moments of $Y$ are estimated. Simplified estimators are available when $b$ has a simple form. For an illustration we consider, in a second example, estimation of the second moment $EY^2 = Er_\vartheta(X)^2 + \sigma^2$. The fully imputed estimator is

$$\frac{1}{n}\sum_{i=1}^n \frac{\sum_{j=1}^n \hat w_j Z_j \{r_{\hat\vartheta}(X_i) + \hat\varepsilon_j\}^2}{\sum_{j=1}^n Z_j} = \frac{1}{n}\sum_{i=1}^n r_{\hat\vartheta}(X_i)^2 + \frac{\sum_{j=1}^n \hat w_j Z_j \hat\varepsilon_j^2}{\sum_{j=1}^n Z_j}.$$

The mean square errors for the fully imputed and the partially imputed estimator (with and without parameter estimation) are given in Table 2.

Consider the lower section on nonlinear regression first. We see that, as expected, the fully imputed estimator outperforms the partially imputed estimator, and that, in part (a) with missing responses, both estimators are far better than the simple estimator in the last column. Using an estimator $\hat{\vartheta}$ for $\vartheta$, or the true value $\vartheta = 2$, does not have much impact on the mean squared error here. The upper half of Table 2 on linear regression, however, shows a different picture: although the mean squared error of the fully imputed and the partially imputed based on the *true* $\vartheta$ are considerably different (which is what we would expect) the values of the estimators based on the ordinary least squares parameter estimator $\hat{\vartheta}$ suggest that the two approaches are asymptotically equivalent. For the extreme case (b) where $\pi(\cdot) = 1$ this would mean that the fully imputed estimator $n^{-1} \sum_{i=1}^{n} r_{\hat{\vartheta}}(X_i)^2 + n^{-1} \sum_{i=1}^{n} \hat{w}_i \hat{\varepsilon}_i^2$ and the empirical estimator $n^{-1} \sum_{i=1}^{n} Y_i^2$ are asymptotically equivalent. This may be surprising but, in fact, it is easy to see that this is exactly what is happening: we consider the special example of linear regression with normal errors and the ordinary least squares estimator $\hat{\vartheta} = \sum_{i=1}^{n} X_i Y_i / \sum_{i=1}^{n} X_i^2$. Rewriting the

TABLE 2
*Simulated mean squared errors of estimators of $EY^2$*

| $\pi(X)$ | $n$ | $\widehat{FI}$ | FI | $\widehat{PI}$ | PI | N |
|---|---|---|---|---|---|---|
| Linear regression: $r_\vartheta(X) = \vartheta X$  ($\vartheta = 2$) | | | | | | |
| $1/(1+e^{-X})$ | 50 | 0.312670 | 0.116360 | 0.310263 | 0.161374 | 0.528146 |
| | 100 | 0.158512 | 0.055343 | 0.157402 | 0.079863 | 0.267601 |
| | 1000 | 0.016215 | 0.005470 | 0.016189 | 0.008113 | 0.027298 |
| 1 | 50 | 0.174683 | 0.070048 | 0.173817 | 0.173817 | 0.173817 |
| | 100 | 0.088960 | 0.034685 | 0.088455 | 0.088455 | 0.088455 |
| | 1000 | 0.008630 | 0.003359 | 0.008623 | 0.008623 | 0.008623 |
| Nonlinear regression: $r_\vartheta(X) = \cos(\vartheta X)$  ($\vartheta = 2$) | | | | | | |
| $1/(1+e^{-X})$ | 50 | 0.086350 | 0.087286 | 0.092361 | 0.093401 | 0.176124 |
| | 100 | 0.042671 | 0.042747 | 0.047054 | 0.047219 | 0.092478 |
| | 1000 | 0.004260 | 0.004179 | 0.005032 | 0.004961 | 0.010153 |
| 1 | 50 | 0.043774 | 0.043873 | 0.066100 | 0.066100 | 0.066100 |
| | 100 | 0.021578 | 0.021574 | 0.035573 | 0.035573 | 0.035573 |
| | 1000 | 0.002159 | 0.002116 | 0.003713 | 0.003713 | 0.003713 |

*Notes.* Here we study estimation of $EY^2$. The first two columns refer to the fully imputed estimator with the ordinary least squares estimator $\hat{\vartheta}$ plugged in ($\widehat{FI}$) and to its version using $\vartheta = 2$ (FI). In the next two columns we consider the partially imputed estimator based on $\hat{\vartheta}$ ($\widehat{PI}$) and $\vartheta = 2$ (PI). In the last column the mean squared errors of $n^{-1} \sum_{i=1}^{n} Z_i Y_i / \pi(X_i)$ (N) are listed.

empirical estimator gives $n^{-1} \sum_{i=1}^n Y_i^2 = n^{-1} \sum_{i=1}^n r_{\hat{\vartheta}}(X_i)^2 + n^{-1} \sum_{i=1}^n \hat{\varepsilon}_i^2 + n^{-1} 2\hat{\vartheta} \sum_{i=1}^n \hat{\varepsilon}_i X_i$. The last term cancels for the least squares estimator $\hat{\vartheta}$ so that $n^{-1} \sum_{i=1}^n Y_i^2 = n^{-1} \sum_{i=1}^n r_{\hat{\vartheta}}(X_i)^2 + n^{-1} \sum_{i=1}^n \hat{\varepsilon}_i^2$. Finally, by our results from Section 3, the estimators $n^{-1} \sum_{i=1}^n \hat{w}_i \hat{\varepsilon}_i^2$ and $n^{-1} \sum_{i=1}^n \hat{\varepsilon}_i^2$ of the error variance $\sigma^2$ are asymptotically equivalent.

In the next example we restrict our attention to linear regression, $r_{\vartheta}(X) = \vartheta X$, $\vartheta = 2$, and consider estimation of a more complicated expectation, namely of $Eh(X, Y) = E(Xe^{XY})$. In contrast to the previous examples the (weighted) fully imputed estimator cannot be reduced. The mean squared errors of this estimator and of the partially imputed estimator are given in Table 3. For each estimator we study the two cases with and without parameter estimation. Again we observe that the performance of the estimators is not much affected by the plug-in parameter estimator. Comparing the fully and the partially imputed estimators we see that the fully imputed estimator clearly outperforms the partially imputed estimator. In addition we also calculate the simulated mean squared error of the unweighted (inefficient) version of our fully imputed estimator. The performance of this estimator turns out to lie between the fully and the partially imputed one. In particular, the simulations in section (b), where all data are observed and where the partially imputed estimator equals the empirical estimator, confirm our theoretical observation that incorporating the information about the location of the errors, for example in the form of weights as done in this article, is important.

In order to study the behavior of the fully imputed estimator for multi-dimensional $\vartheta$ we again studied estimation of $E(Xe^{XY})$. For our simulations

TABLE 3
*Simulated mean squared errors of estimators of $E\{X \exp(XY)\}$ in linear regression*

| $\pi(X)$ | $n$ | $\widehat{\mathbf{FI}}$ | $\mathbf{FI}$ | $\widehat{\mathbf{U}}$ | $\widehat{\mathbf{PI}}$ | $\mathbf{PI}$ |
|---|---|---|---|---|---|---|
| $1/(1 + e^{-X})$ | 50 | 0.32563 | 0.29024 | 0.36187 | 0.48164 | 0.47769 |
| | 100 | 0.15017 | 0.14085 | 0.18147 | 0.24192 | 0.24698 |
| | 1000 | 0.01384 | 0.0137 | 0.01992 | 0.02577 | 0.02703 |
| 1 | 50 | 0.28988 | 0.27262 | 0.32220 | 0.58566 | 0.58566 |
| | 100 | 0.13804 | 0.13413 | 0.16520 | 0.29948 | 0.29948 |
| | 1000 | 0.01332 | 0.01329 | 0.01663 | 0.02997 | 0.02997 |

*Notes.* We consider estimation of $Eh(X, Y) = E(Xe^{XY})$ in the linear regression model $r_{\vartheta}(X) = \vartheta X$, $\vartheta = 2$. The first two columns give the mean squared errors of the fully imputed estimator with the least squares estimator $\hat{\vartheta}$ plugged in ($\widehat{\text{FI}}$), and its version using $\vartheta = 2$ (FI). The third column contains the mean squared errors of the unweighted version $\widehat{\text{U}}$ of $\widehat{\text{FI}}$. The last two columns refer to the partially imputed estimator using $\hat{\vartheta}$ ($\widehat{\text{PI}}$) and $\vartheta = 2$ (PI). Note that if $\pi(X) = 1$ then the partially imputed estimator again equals the empirical estimator, $\text{PI} = \widehat{\text{PI}} = n^{-1} \sum_{i=1}^n X_i \exp(X_i Y_i)$.

Table 4

*Simulated mean squared errors of estimators of $E\{X\exp(XY)\}$ with $\vartheta \in \mathbb{R}^p$ $(p=2,3)$*

| $r_\vartheta(X)$ | $\widehat{\text{FI}}$ | **FI** | $\widehat{\text{U}}$ | $\widehat{\text{PI}}$ | **PI** |
|---|---|---|---|---|---|
| $\vartheta_1 X + \vartheta_2 U$ | 0.2465 | 0.2272 | 0.2855 | 0.3965 | 0.4018 |
| $\vartheta_0 + \vartheta_1 X + \vartheta_2 U$ | 0.3048 | 0.2272 | 0.3048 | 0.4259 | 0.4018 |
| $\vartheta_1 X + \vartheta_2 U + \vartheta_3 V^2$ | 0.4367 | 0.3750 | 0.4434 | 0.5696 | 0.5760 |

*Notes.* The three rows refer to two regression functions with different parametrizations. We have $\vartheta_0 = 0$, $\vartheta_1 = 2$, $\vartheta_2 = -1$ and $\vartheta_3 = 0.5$, $n = 100$, $\pi(X) = 1/(1 + e^{-X})$. The covariates $X, U$ and $V$ are independent from a uniform distribution on $(-1, 1)$. The parameters are estimated using least squares. The notation is explained in Table 3.

Table 5

*Simulated mean squared errors of estimators of $Eh(X,Y)$ with $\hat{\vartheta}$ inefficient*

| | | **$EY$** $r_\vartheta(X) = \cos(\vartheta X)$ | | **$EY$** $r_\vartheta(X) = \vartheta X$ | | **$E(Xe^{XY})$** $r_\vartheta(X) = \vartheta X$ | |
|---|---|---|---|---|---|---|---|
| $\pi(X)$ | $n$ | $\widehat{\text{FI}}$ | $\widehat{\text{PI}}$ | $\widehat{\text{FI}}$ | $\widehat{\text{PI}}$ | $\widehat{\text{FI}}$ | $\widehat{\text{PI}}$ |
| $1/(1 + e^{-X})$ | 50 | 0.03124 | 0.03545 | 0.02742 | 0.03868 | 0.50275 | 0.72944 |
| | 100 | 0.01841 | 0.02057 | 0.01375 | 0.01938 | 0.24148 | 0.48759 |
| 1 | 50 | 0.02000 | 0.02864 | 0.02689 | 0.05181 | 0.41476 | 0.79949 |
| | 100 | 0.01016 | 0.01448 | 0.01359 | 0.02589 | 0.25796 | 0.63103 |

*Notes.* We compare fully and the partially imputed estimators of $EY$ and $E(Xe^{XY})$, keeping the previous notation. Again, $\hat{\vartheta}$ is the least squares estimator, but now the errors are from a $t$-distribution with 10 degrees of freedom.

we restricted our attention to missing data and on samples of size $n = 100$, and considered three different regression models which are given in Table 4. Note that the second regression function, $\vartheta_0 + \vartheta_1 X + \vartheta_2 U$ with $\vartheta_0 = 0$, $\vartheta_1 = 2$ and $\vartheta_2 = -1$, equals the first one, namely $2X - U$, but it involves a three-dimensional parameter. As expected, the increase of dimension impairs the performance of the fully imputed (weighted and unweighted) and of the partially imputed estimator. Note that the weighted and unweighted fully imputed estimator ($\widehat{\text{FI}}$ and $\widehat{\text{U}}$) in the second regression model are the same: we consider the least squares estimator in a regression model with an intercept term $\vartheta_0$. In this model the least squares estimator solves, by construction, $\sum_{j=1}^{n} Z_j \hat{\varepsilon}_j = 0$ (which implies that all weights $\hat{w}_j$ equal one). Again we observe that the fully imputed estimator consistently outperforms the partially imputed estimator.

We conclude this section with a small simulation study to examine the behavior of the fully imputed estimator when $\hat{\vartheta}$ is inefficient. The simplest setting is to choose the ordinary least squares estimator, as we did before, but in a model with non-normal errors. In Table 5 we consider estimation of the mean response and of $E(Xe^{XY})$ for linear and nonlinear regression, and for errors from a $t$-distribution. The results are similar to the previous ones: again the fully imputed estimator performs best, though not as well as if the errors are, in fact, from a normal distribution (cf. Tables 1–3). Simulations with a logistic error density turned out similarly, confirming the better performance of the imputation method. At least in these examples, with moderate sample sizes $n = 50$ and $n = 100$, the construction of $\vartheta$ does not seem to be as important as the choice between the full and the partial imputation approaches.

6.3. *Confidence intervals.* By Theorem 5.2 the fully imputed weighted estimator $n^{-1} \sum_{i=1}^{n} \hat{\chi}_w(X_i, \hat{\vartheta})$ is asymptotically normally distributed, with asymptotic variance $\sigma_{\mathrm{FI}}^2 = E\chi^2(X, \vartheta) + (EZ)^{-1}E\bar{h}^2(\varepsilon) - \{1 + (EZ)^{-1}\}E^2h(X, Y) - E^2\{\varepsilon\bar{h}(\varepsilon)\}/(\sigma^2 EZ) + D_w^\top E(Z\zeta\zeta^\top)^{-1}D_w$ (see Theorem 5.2 for the notation). An asymptotic confidence interval for $Eh(X, Y)$ with confidence level $1 - \alpha$ is

$$\left( \frac{1}{n}\sum_{i=1}^{n} \hat{\chi}_w(X_i, \hat{\vartheta}) - z_{\alpha/2}\sqrt{\frac{\hat{\sigma}_{\mathrm{FI}}^2}{n}}, \frac{1}{n}\sum_{i=1}^{n} \hat{\chi}_w(X_i, \hat{\vartheta}) + z_{\alpha/2}\sqrt{\frac{\hat{\sigma}_{\mathrm{FI}}^2}{n}} \right),$$

where $z_{\alpha/2}$ denotes the upper $\alpha/2$-quantile of the standard normal distribution, and where $\hat{\sigma}_{\mathrm{FI}}^2$ is a consistent estimator of $\sigma_{\mathrm{FI}}^2$. Consider, for example, estimation of $EY$ with $r_\vartheta(X)$ depending on a scalar parameter $\vartheta$, which covers our previous simple examples $r_\vartheta(X) = \vartheta X$ and $r_\vartheta(X) = \cos(\vartheta X)$. Here the confidence interval is $n^{-1}\sum_{i=1}^{n} r_{\hat{\vartheta}}(X_i) \pm z_{\alpha/2}(\hat{\sigma}_{\mathrm{FI}}^2/n)^{1/2}$. The asymptotic variance of $n^{-1}\sum_{i=1}^{n} r_{\hat{\vartheta}}(X_i)$ is

$$\sigma_{\mathrm{FI}}^2 = \mathrm{Var}\, r_\vartheta(X) + \frac{E^2\{\dot{r}_\vartheta(X)\}}{EZ\,\mathrm{Var}\{r_\vartheta(X)|Z = 1\}E\{\ell^2(\varepsilon)\}}.$$

The expectations in the formula can be estimated by empirical methods, with a consistent estimator $\hat{\vartheta}$ for the parameter $\vartheta$ plugged in. Consider, for example, $\mathrm{Var}\{r_\vartheta(X)|Z = 1\} = E\{r_\vartheta(X)^2|Z = 1\} - E^2\{r_\vartheta(X)|Z = 1\}$. The first expectation is estimated by $(\sum_{i=1}^{n} Z_i)^{-1}\sum_{i=1}^{n} Z_i\{r_{\hat{\vartheta}}(X_i)\}^2$, and analogously the second one.

In order to confirm the theoretical results we also performed some simulation studies, generating confidence intervals for the above examples with the described estimation method. As expected, for $\alpha = 0.05$ we obtained the desired coverage probability 0.95.

## APPENDIX

LEMMA A.1. *Suppose that Assumption 1 is satisfied. Then, for a $\sqrt{n}$ consistent estimator $\hat{\vartheta}$ of $\vartheta$, the statements* (3.2)–(3.4) *hold.*

PROOF.   In order to prove (3.2)–(3.4) we first show

(A.1)    $\max_{1 \leq i \leq n} |Z_i \hat{\varepsilon}_i - Z_i \varepsilon_i| = o_p(1),$

(A.2)    $\sum_{i=1}^{n} Z_i (\hat{\varepsilon}_i - \varepsilon_i^*)^2 = o_p(1) \qquad \text{with } \varepsilon_i^* = \varepsilon_i - \dot{r}_\vartheta(X_i)^\top (\hat{\vartheta} - \vartheta).$

Result (A.2) immediately follows from the $\sqrt{n}$ consistency of $\hat{\vartheta}$ and the stochastic differentiability of $r_\vartheta$ [implication (2.1) of Assumption 1]:

$$\sum_{i=1}^{n} Z_i (\hat{\varepsilon}_i - \varepsilon_i^*)^2 = \sum_{i=1}^{n} Z_i [\hat{\varepsilon}_i - \{\varepsilon_i - \dot{r}_\vartheta(X_i)^\top (\hat{\vartheta} - \vartheta)\}]^2$$

$$\leq \sum_{i=1}^{n} \{r_{\hat{\vartheta}}(X_i) - r_\vartheta(X_i) - \dot{r}_\vartheta(X_i)^\top (\hat{\vartheta} - \vartheta)\}^2 = o_p(1).$$

This gives $\max_{1 \leq i \leq n} |Z_i (\hat{\varepsilon}_i - \varepsilon_i^*)| = o_p(1)$. In order to establish (A.1) it therefore suffices to show $\max_{1 \leq i \leq n} |Z_i (\varepsilon_i^* - \varepsilon_i)| = o_p(1)$. We have

$$\max_{1 \leq i \leq n} |Z_i (\varepsilon_i^* - \varepsilon_i)| \leq \max_{1 \leq i \leq n} |\varepsilon_i^* - \varepsilon_i| \leq |\hat{\vartheta} - \vartheta| \cdot \max_{1 \leq i \leq n} |\dot{r}_\vartheta(X_i)|.$$

Since $\hat{\vartheta}$ is $\sqrt{n}$ consistent we only need $n^{-1/2} \max_{1 \leq i \leq n} |\dot{r}_\vartheta(X_i)| = o_p(1)$. But this holds by Owen (2001), Lemma 11.2, since the variables $|\dot{r}_\vartheta(X_i)|$, $i = 1, \ldots, n$, are i.i.d. and, by Assumption 1, have finite second moments. This shows $\max_{1 \leq i \leq n} |Z_i (\varepsilon_i^* - \varepsilon_i)| = o_p(1)$.

Equation (3.2), $\max_{1 \leq i \leq n} |Z_i \hat{\varepsilon}_i| = o_p(n^{1/2})$, can be seen as follows: we can bound $\max_{1 \leq i \leq n} |Z_i \hat{\varepsilon}_i|$ by $\max_{1 \leq i \leq n} |Z_i \hat{\varepsilon}_i - Z_i \varepsilon_i| + \max_{1 \leq i \leq n} |Z_i \varepsilon_i|$. The first term is $o_p(1)$ by (A.1) and the second term is $o_p(n^{1/2})$ by Owen's Lemma 11.2 since the $Z_i \varepsilon_i$ are i.i.d. with finite variance. We now show (3.3), that is,

$$\frac{1}{n} \sum_{i=1}^{n} Z_i \hat{\varepsilon}_i = \frac{1}{n} \sum_{i=1}^{n} Z_i \varepsilon_i - E Z E\{\dot{r}_\vartheta(X)|Z = 1\}^\top (\hat{\vartheta} - \vartheta) + o_p(n^{-1/2}).$$

In view of (A.2), $n^{-1} \sum_{i=1}^{n} Z_i \hat{\varepsilon}_i = n^{-1} \sum_{i=1}^{n} Z_i \varepsilon_i^* + o_p(n^{-1/2})$. By the law of large numbers we obtain

$$\frac{1}{n} \sum_{i=1}^{n} Z_i \varepsilon_i^* = \frac{1}{n} \sum_{i=1}^{n} Z_i \varepsilon_i - \frac{1}{n} \sum_{i=1}^{n} Z_i \ddot{r}_\vartheta(X_i)^\top (\hat{\vartheta} - \vartheta)$$

$$= \frac{1}{n} \sum_{i=1}^{n} Z_i \varepsilon_i - E\{Z \dot{r}_\vartheta(X)\}^\top (\hat{\vartheta} - \vartheta) + o_p(n^{-1/2}).$$

Since $E\{Z \dot{r}_\vartheta(X)\} = E Z E\{\dot{r}_\vartheta(X)|Z = 1\}$ we have established (3.3).

Our last auxiliary result to prove is (3.4),

$$\frac{1}{n} \sum_{i=1}^{n} Z_i \hat{\varepsilon}_i^2 = \frac{1}{n} \sum_{i=1}^{n} Z_i \varepsilon_i^2 + o_p(1) = E Z \sigma^2 + o_p(1).$$

The second equality is just a consequence of the law of large numbers. To see that the first equation holds consider

$$\frac{1}{n}\sum_{i=1}^{n}Z_i\hat{\varepsilon}_i^2 - \frac{1}{n}\sum_{i=1}^{n}Z_i\varepsilon_i^2 = \frac{1}{n}\sum_{i=1}^{n}Z_i(\hat{\varepsilon}_i-\varepsilon_i)^2 + 2\frac{1}{n}\sum_{i=1}^{n}Z_i(\hat{\varepsilon}_i-\varepsilon_i)\varepsilon_i.$$

The second term on the right-hand side is $o_p(1)$ by (A.1). To show that the first expression is $o_p(1)$ it suffices, in view of (A.2), to consider

$$\frac{1}{n}\sum_{i=1}^{n}Z_i(\varepsilon_i^*-\varepsilon_i)^2 = \frac{1}{n}\sum_{i=1}^{n}Z_i\{\dot{r}_{\vartheta}(X_i)^{\top}(\hat{\vartheta}-\vartheta)\}^2.$$

This term is $o_p(1)$ since $\hat{\vartheta}$ is $\sqrt{n}$ consistent and since $\dot{r}_{\vartheta}(X)$ is in $L_2(P)$.   $\square$

**Acknowledgments.**   Many thanks to Anton Schick for important suggestions on specifying and constructing efficient parameter estimators in Section 5, to Wolfgang Wefelmeyer for valuable discussions and advice and to Raymond Carroll for constructive criticism of an earlier draft. Thanks also to two referees for helpful comments.

DEPARTMENT OF STATISTICS
TEXAS A&M UNIVERSITY
COLLEGE STATION, TEXAS 77843-3143
USA
E-MAIL: uschi@stat.tamu.edu
URL: http://www.stat.tamu.edu/~uschi/